\documentclass[leqno,11pt]{article}

\usepackage{amsfonts,amssymb,amsgen,theorem,stmaryrd, amsopn}
\usepackage[leqno]{amsmath} 
\usepackage[T1]{fontenc}
\usepackage[latin1]{inputenc} 
\usepackage[french,english]{babel}
\usepackage{graphicx}
\usepackage{graphpap}

\newtheorem{theoreme}{Theorem} \theorembodyfont{\sl}
\newtheorem{proposition}{Proposition} 
\newtheorem{lemme}[proposition]{Lemma}

\newtheorem{rem}[proposition]{Remark} 

\theoremheaderfont{\sc}
\numberwithin{equation}{section}
\numberwithin{theoreme}{section}
\numberwithin{proposition}{section}

\newcommand{\eps}{\varepsilon} 

\renewcommand\Re{\mathrm{Re}\,} \renewcommand\Im{\mathrm{Im}\,}
\newcommand\R{{\mathbb R}}

 \renewcommand\S{{\mathbb S}}
\newcommand\supp{{\mathrm{supp}\hskip 0.05cm}}

\global
 \global
 \theoremstyle{break} \theorembodyfont{\it}

\renewcommand{\Im}{  \text{Im}   }
\renewcommand{\Re}{  \text{Re}   }

 \def\imp{\Longrightarrow}
 \def\tq{\,\,\mbox{s.t.}\,\,}

 \def\cdotv{\raise 2pt\hbox{,}}

\begin{document}

 \title{\vskip -1cm Bilinear virial
   identities and applications}

\author{Fabrice Planchon  \footnote{The first author was partially supported by A.N.R. grant ONDE NON
   LIN}\\Laboratoire Analyse, G\'eom\'etrie
   \& Applications, UMR 7539 du CNRS,\\ Institut Galil\'ee, Universit\'e Paris
   13,\\ 99 avenue J.B. Cl\'ement, F-93430 Villetaneuse\\fab@math.univ-paris13.fr\\
\and
Luis Vega  \footnote{The second author was partially supported by grant MTM2007-62186}
 \\Universidad del Pais Vasco - Euskal Herriko Unibertsitatea,\\
 Departamento de Matematicas Aptdo. 644,\\
  48080 Bilbao\\
luis.vega@ehu.es}

\date{}
 \maketitle
 \begin{abstract}
We prove bilinear virial identities for the nonlinear Schr\"odinger equation,
which are extensions of the Morawetz interaction inequalities. We
recover and extend known bilinear improvements to Strichartz
inequalities and provide applications to various nonlinear problems,
most notably on domains with boundaries.
\\
\selectlanguage{french}
\begin{center}
  Résumé
\end{center}
On démontre des identités de type viriel bilinéaire pour l'équation de
Schr\"odinger nonlinéaire, qui peuvent être vues comme des extensions
des inéqualités d'interaction de Morawetz. Ceci permet de retrouver et
d'étendre des rafinements bilinéaires des inéqualités de Strichartz,
et nous donnons également des applications à plusieurs problèmes
nonlinéaires, notamment sur les domaines à bord.
\selectlanguage{english}
\end{abstract}
 \par \noindent
\section{Introduction}
Dispersive estimates are known to be an essential tool in dealing with
low regularity well-posedness issues for the nonlinear Schr\"odinger
equation. Among the most useful ones are Strichartz inequalities: starting
with \cite{Stri}, they were completed by \cite{GV} and finally by
\cite{KT}. As space-time bounds for solutions to the linear
Schr\"odinger equation in $\R^n$, they are closely related to the Fourier
restriction problem in harmonic analysis
, and as such
heavily rely on the use of Fourier transform techniques. Extensions of these
inequalities to more complicated geometrical settings have been the
subject of intense research over the last decade, to the point where
quoting all possible references would fill this page. It should be
noted that these works are based on appropriate refinements of the
$\R^n$ case, through Fourier Integral Operator, FBI, wave packet or any
appropriate microlocal generalizations of Fourier analysis (for a
notable exception using vector field methods, see
\cite{Salort}). On the other hand, one has virial type identities, of
which the Morawetz identity (proved by Lin-Strauss \cite{LS}) is
perhaps the most well-known: such identities have two key features,
they are obtained by integration by parts and they usually apply to
the nonlinear equation. We remark that the local smoothing effect, which came much later and was first observed in the flat case (see \cite{CS}, \cite{Sj}, \cite{Luis}), may be
seen as part of this category as well, though proofs usually require a
sophisticated ``integration by parts'' involving pseudo-differential
operators or resolvent methods. A new kind of inequality was introduced in
\cite{CKSTTmora}, the Morawetz interaction inequality, which seemed to
have the benefit of both worlds: one may recover a specific,
non-sharp Strichartz estimate and it also applies to the nonlinear equation
(providing an essential tool to solve the $H^1$-critical defocusing
NLS in $3D$, \cite{CKSTTannals}). Subsequent developments include a
curved space version (\cite{HTW}) and a quartic interaction inequality
for NLS on $\R$ (\cite{CHV}).

In the present work, we explore a different direction, which builds
upon the understanding of the local smoothing effect and its
fundamentally $1D$ nature. This naturally leads to a new set of
identities with several interesting consequences:
\begin{itemize}
\item in $1D$, one recovers, by a simple argument, an identity of
  \cite{OzTsu}, which implies the Fefferman-Stein inequality in its
  bilinear version; from there the (almost) full set of
  Strichartz/maximal function estimates may be derived. More
  importantly, we get a nonlinear identity.
\item In $2D$ and higher, one obtains an $L^2_{t,x}$-based estimate
  for the charge density. (This would correspond, w.r.t. scaling, to a sharp Strichartz estimate in $2D$). More interestingly, one may derive from our result Bourgain's bilinear improvement
  (\cite{Bourgain98}).
\item All our identities apply to nonlinear equations, and have
  bilinear versions.
\item Nothing but integration by parts is used in the proof: as such,
  these estimates extend to domains, provided one may control the
  boundary terms; in the case of Dirichlet boundary conditions, such
  control is provided by local smoothing.
\item As an application to exterior domains, we improve the
  well-posedness theory to $H^1$-subcritical (subquintic) nonlinearities for $n=3$.
\item Applications to scattering problems are straightforward, and this
  extends to 3D exterior domains, where no results were
  available to our knowledge and where we obtain scattering in the
  energy class for the defocusing cubic equation. 
\end{itemize}
While presenting this work at Oberwolfach, we learned that similar
results (namely a priori bound \eqref{eq:borneapriori}) have been obtained
simultaneously and independently by J. Colliander, M. Grillakis and
N. Tzirakis, see \cite{CGT} and \cite{CGTmora}), through a different derivation.
\vskip 0.2cm
\par\noindent
{\bf Acknowledgments:} we thank N. Burq for various enlightenments
about the Schr\"odinger equation on exterior domains, as well as the
referee for helpful comments and suggestions which greatly improved
the presentation.
\section{Main results}
\subsection{The Schr\"odinger equation in $\R^n$}
Let $n\geq 1$, $p\in \R$, $p\geq 1$, $\eps\in \{-1,0,1\}$, and $u$
is a solution to 
\begin{equation}
\label{equ}
  i\partial_t u+ \Delta u=\eps |u|^{p-1} u, \text{ with } u_{|t=0}=u_0.
\end{equation}
We will also need $v$, solution to
\begin{equation}
\label{eqv}
  i\partial_t v+ \Delta v=\eps |v|^{p-1} v, \text{ with } v_{|t=0}=v_0.
\end{equation}
Let us define several quantities which will play a key role:
for $n>1$ and given a function $f$, its Radon transform is
\begin{equation}
  R(f)(s,\omega)=\int_{x\cdot \omega=s} f \,d\mu_{s,\omega},
\end{equation}
where $\mu_{s,\omega}$ is the induced measure on the hyperplane
$x\cdot\omega=s$. We set
\begin{equation}
  \label{eq:virial1}
  I_\omega (\eps,u,v)=\int_{x\cdot \omega> y\cdot \omega}
  (x\cdot \omega- y\cdot \omega) |u|^2(x) |v|^2(y)\,dxdy.
\end{equation}
Remark that a simple computation leads to
\begin{equation}
\partial_t I_\omega =i \left (\int_{x\cdot \omega >y\cdot
    \omega}\omega\cdot\left [   ( u \nabla_x \bar u-\bar u
  \nabla_x u)(x)|v(y)|^2-  ( v \nabla_y \bar v-\bar v
  \nabla_y v)(y)|u(x)|^2\right ] \,dy\,dx\right).
\end{equation}
We may now state our first result.
\begin{theoreme}
  \label{t1}
Let $\omega \in \R^n$, $n>1$, with $|\omega|=1$, $u$ solution to
\eqref{equ}. Then, with $x=x^\perp+s\omega$
\begin{multline}
\label{nonlinRad}
\int_s |\partial_s (R(|u|^2))(s,\omega)|^2\,ds+\eps \frac{p-1}{p+1} \int_s R(|u|^2)
R(|u|^{p+1})\,ds\\
{}+\int_s \int_{x^\perp\cdot \omega=0}\int_{y^\perp \cdot \omega=0} |u(x^\perp+s\omega)\partial_{s}
u(y^\perp+s \omega)-u(y^\perp+s\omega)\partial_{s}
u(x^\perp+s \omega)|^2 \,dx^\perp dy^\perp  ds \\
= \frac 1 4 \partial^2_t I_\omega(\eps,u,u) 
\end{multline}
In other words, $I_\omega(\eps,u,u)$ is a convex function in time.
\end{theoreme}
In the specific $1D$ case, one has actually the following identity.
\begin{theoreme}
  \label{t2}
Let $n=1$, $u,v$ two solutions to \eqref{equ}, \eqref{eqv}, then
\begin{equation}
  \label{eq:nlid1D}
  4\int_x |\partial_x (u \bar v)|^2 \, dx+2 \eps \frac {p-1}{p+1} \int_x
  |u|^2 |v|^{p+1}+|v|^2 |u|^{p+1} \,dx = \partial^2_t I(\epsilon,u,v).
\end{equation}
\end{theoreme}
\begin{rem}
Up to a doubling factor, $I_\omega$ may be
  recast as a Morawetz interaction functional (as introduced in \cite{CKSTTmora}),
$$
\int \rho(x-y) |u|^2(x) |v|^2(y)\, dxdy,
$$
with $\rho(x-y)=|x\cdot \omega -y\cdot \omega|$. Hence we have
replaced the physical distance $|x-y|$ (which was the default choice
in \cite{CKSTTmora} and subsequent works) by its projection over a
specified direction $\omega$. We chose our definition of $I_\omega$ as to
emphasize trace terms which will later appear in the proof. In fact,
we were led to $I_\omega$ by considering variations on the local
smoothing, and we will come back to this point in section \ref{localsmoothingsection}.
\end{rem}
In order to turn these bounds into useful nonlinear control, we use
\begin{proposition}
\label{propH12}
  Let $\omega$ be fixed, then
  \begin{equation}
    \label{eq:moment}
    |\partial_t I_\omega|\leq \|u\|^2_{L^2_x}\|v\|^2_{\dot H^\frac 1
     2}+\|v\|^2_{L^2_x}\|u\|^2_{\dot H^\frac 1 2}.
  \end{equation}
As a consequence, when $\eps=1$ (defocusing equation), we have an
priori bound,
\begin{multline}
  \label{eq:borneapriori}
  \int_{\R} \int_{\R^n} | |\nabla|^\frac {3-n} 2 (|u|^2)|^2 \,dx dt+  \int_{\R} \int_{\R^n} | |\nabla|^\frac {1-n} 2 (|u|^{\frac{p+3}2})|^2 \,dx dt  \lesssim sup_{t\in \R} \|u\|^2_{L^2_x}\|u\|^2_{\dot H^\frac 1
     2}.
\end{multline}
\end{proposition}
  \begin{rem}
    The right-hand side of \eqref{eq:moment} is very clearly not invariant by
    galilean transforms. The left-hand side, however, is.
  \end{rem}
  \begin{rem}
    The a priori estimate \eqref{eq:borneapriori} was obtained 
simultaneously and independently by J. Colliander, M. Grillakis and
N. Tzirakis (\cite{CGT} and \cite{CGTmora}), through a direct
derivation with the weight $\rho(x)=|x|$ but with a new commutator
argument involving $[x,\sqrt{-\Delta}^{-(n-1)}]$ and the local
conservation laws for mass and momentum densities, overcoming
the restriction to dimensions $n\geq 3$ from \cite{CKSTTmora}.
  \end{rem}
We now state a more general result: let
\begin{equation}
  \label{eq:imf}
  I_\rho(u,v)=\int \rho(x-y) |u|^2(x) |v|^2(y)\,dxdy.
\end{equation}
Then
\begin{theoreme}
\label{t3}
Let $\rho$ be a weight function such that its Hessian $H_\rho$ is positive; let
\begin{equation}
  F(u,v)(x,y)=\bar v(y)\nabla_x
 u(x)+u(x)\nabla_y \bar v(y) \text{ and } G(u,v)(x,y)= v(y)\nabla_x
 u(x)-u(x)\nabla_y  v(y).
\end{equation}
 We have
  \begin{eqnarray}
\label{nonlinearmauvais}
 \partial_t^2 I_\rho & = & 4\int H_\rho(x-y)(F(u,v)(x,y),\overline F(u,v)(x,y)) \,dxdy\nonumber\\
 & & {}+ \eps \frac{p-1}{p+1} \int |v|^2(y) (\Delta_x \rho)(x-y)  |u|^{p+1} (x)
  \,dxdy\\
& & {}+ \eps \frac{p-1}{p+1} \int |u|^2(x) (\Delta_x \rho)(x-y)  |v|^{p+1} (y) \,dxdy.\nonumber
\end{eqnarray}
Moreover, we may rewrite
\begin{multline}
\label{nonlinearrho}
\int H_\rho(x-y)(F(u,v)(x,y),\overline F(u,v)(x,y))\,dxdy  = \\
\int
H_\rho(x-y)(G(u,v)(x,y),\overline G(u,v)(x,y))\,dxdy+ \int \Delta
\rho(x-y)\nabla_x(|u|^2(x))\cdot \nabla_y(|v|^2(y)))\,dxdy.
\end{multline}
\end{theoreme}
\begin{rem} Notice that if we make $u=v$ in \eqref{nonlinearrho} and
  assume that the Fourier transform of  $\Delta \rho$ is positive, we can bound each of the two terms in
  the r.h.s. in terms of the l.h.s.
\end{rem}
The above remark used in the particular case $\rho(z)=|z\cdot\omega|$
gives us the following Corollary for the linear equation.
\begin{theoreme}
\label{bornebilineaire}
 Given $\omega$ a unitary vector in $\R^n$,  $n>1$ and $u,v$ any two
 solutions to \eqref{equ} and \eqref{eqv} with $\epsilon=0$ (linear equation) we have
\begin{equation}
 \int_t \int_s |\partial_s (R(u\bar v))(s,\omega)|^2\,ds dt \lesssim I_\omega(\hat u_0,\hat u_0) +I_\omega(\hat v_0,\hat v_0)+I_\omega(\hat u_0,\hat v_0).
\end{equation}
with $I_\omega$ as given in \eqref{eq:virial1}.
\end{theoreme}
\begin{rem}
  We will see that this bilinear estimate implies Bourgain's bilinear
  refinement of Strichartz estimate from \cite{Bourgain98}. One may
  notice that \eqref{bornebilineaire} (and all identities involving
  the Radon transform) does not depend on the
  dimension $n$, in sharp contrast with \eqref{eq:borneapriori}, which
  gets worse with $n$ large.
\end{rem}
\subsection{The Schr\"odinger equation on a domain $\Omega$}
Let $n\geq 1$, $p\in \R$, $p\geq 1$, $\eps\in \{-1,0,1\}$,
$\Omega\subset \R^n$ with a smooth boundary $\partial\Omega$, and $u$
is now the solution to 
\begin{equation}
\label{equd}
  i\partial_t u+\Delta u=\eps |u|^{p-1} u, \text{ with }
  u_{|\partial \Omega} =0.
\end{equation}
Denote by
\begin{equation}
  \label{eq:massenergy}
  M(u)=\int_\Omega |u|^2 \,dx \text{ and } E(u)=\frac 1 2 \int_\Omega |\nabla u|^2 \,
dx+\frac 1 {p+1} \int_\Omega |u|^{p+1} \,dx,
\end{equation}
the mass and energy which are conserved quantities: we will use $M$
and $E$ as shorter version of $M(u)$ and $E(u)$. Notice that the Radon transform is still defined,
\begin{equation}
  R(f)(s,\omega)=\int_{x\cdot \omega=s \cap \Omega} f \,d\mu_{s,\omega}.
\end{equation}
We set
\begin{equation}
\label{eq:virial1bis}
  I_\rho =\int_{x,y\in \Omega}
  \rho(x-y) |u|^2(x) |u|^2(y)\,dxdy.
\end{equation}
We may now state our result.
\begin{theoreme}
  \label{t1d}
Let $\omega \in \R^n$, $n>1$, with $|\omega|=1$, and pick $\rho_\omega (z)=|z\cdot\omega|$, $u$ solution to
\eqref{equd}. Then, with $x=x^\perp+s\omega$
\begin{multline}
\label{idt1}
\int_s |\partial_s (R(|u|^2))(s,\omega)|^2\,ds+\eps \frac{p-1}{p+1} \int_s R(|u|^2)
R(|u|^{p+1})\,ds\\
{}+\int_s \int_{x\cdot \omega=s}\int_{y\cdot \omega=s} |u(x^\perp+s\omega)\partial_{s}
u(y^\perp+s \omega)-u(y^\perp+s\omega)\partial_{s}
u(x^\perp+s \omega)|^2 \,dx^\perp dy^\perp  ds \\
 -\int_{x\in
  \partial\Omega, y\in\Omega} |u|^2(y) \partial_n \rho_\omega(x-y)
  |\partial_n u|^2(x) \, dS_x dy = \partial^2_t I_{\rho_\omega}.
\end{multline}
\end{theoreme}
We now illustrate how to obtain useful estimates from Theorem
\ref{t1d} when one has control of the boundary term.
\begin{proposition}
\label{p1d}
  Let $\Omega$ be $\R^n\setminus \Sigma$, where $\Sigma$ is
  star-shaped and $\Sigma\subset\subset K$, $K$ compact. Assume moreover $\eps=0,1$ (linear or defocusing) and
  $n\geq 3$. Then, 
  \begin{equation}
    \label{eq:nlsdnormal}
\int_0^T \int_{x\in \partial\Omega} |\partial_n u|^2 \,dS_x dt
+\int_0^T\int_{K\setminus \Sigma} (|\nabla u|^2+|u|^2) \,dx dt\lesssim
\sup_{t\in [0,T]} \|u\|^2_{\dot H_0^\frac 1 2
  (\Omega)}\lesssim (ME)^\frac 1 2.
  \end{equation}
\end{proposition}
\begin{rem}
In 2D, one may only obtain a local in time estimate for the defocusing
equation, an issue related to the zero mode (or the failure of the Morawetz
estimate in $2D$) Hence,  \eqref{eq:nlsdnormal} will have an
additional term $C(T) \|u_0\|^2_{L^2}$ on the right-hand side. We will
not use such an estimate and therefore skip it.
\end{rem}
As a consequence of Theorem \ref{t1d} and Proposition \ref{p1d}, we
have
\begin{proposition}
  \label{p2d}
  Let $\Omega$ be $\R^n\setminus \Sigma$, where $\Sigma$ is
  star-shaped and $\Sigma\subset\subset K$, $K$ compact, and  $n\geq
  3$. Then, the solution $u$ to the defocusing ($\varepsilon=1$)
  equation\eqref{equd} verifies
  \begin{equation}
    \label{eq:nlcd}
    \||\nabla_x|^{\frac{3-n} 2}(|u|^2)\|_{L^2_{t,x}}\lesssim \sup_t
    \|u\|_{L^2(\Omega)}\|u\|_{\dot H_0^\frac 1 2} \lesssim M^\frac 3 4
    E^\frac 1 4.
  \end{equation}
\end{proposition}
Note that, more generally, for the linear equation, the result of Proposition \ref{p1d} holds
for unbounded domains, assuming one does not have any
trapped rays. In fact, for such domains, the local smoothing estimate
holds (\cite{BurqDuke}), irrespective of the dimension and with an
absolute constant (independent of $T$); a simple integration by part argument (close to the
boundary) yields control of the boundary term. As such, one obtains
\begin{theoreme}
  \label{t2d}
Let $n\geq 2$, $\Omega \subset \R^n$ be a domain where
\eqref{eq:nlsdnormal} holds for the linear equation, and $u$ a solution to
the linear equation \eqref{equd} ($\eps=0$). Then the following
estimate holds: 
\begin{equation}
\||\nabla_x|^{\frac{3-n} 2}(|u|^2)\|_{L^2_{t,x}}\lesssim \|u_0 \|^2_{\dot H_0^\frac 1 4(\Omega)}.
\end{equation}
\end{theoreme}
Now, consider the linear equation on a domain for which local
smoothing does not hold. By tailoring the size of the time interval to
the frequency of the solution, one may obtain an estimate with a $1/4$
loss of regularity.
\begin{theoreme}
  \label{t2dbis}
Let $\Omega \subset \R^n$ be a bounded domain, $u$ solution to
the linear equation \eqref{equd} ($\eps=0$). Then
\begin{equation}
\label{L4borne}
\||\nabla_x|^{\frac{3-n} 2}(|u|^2)\|_{L^2([0,1];L^2(\Omega))}\lesssim
\|u_0 \|^2_{\dot H_0^\frac 1 2(\Omega)}.
\end{equation}
\end{theoreme}
We remark that the boundedness of the domain is in no way essential.
\begin{rem}
  The numerology of \eqref{L4borne} is consistent with the numerology
  of \cite{BuGeTz01-1} on manifolds without boundaries. By contrast,
  estimates from \cite{blair} have an additional $\frac 1 {3p}$ loss, where $p$ is
  the time Lebesgue exponent; our example suggests better estimates
  than the ones which are obtained by interpolation between the $p=2$
  case and the conservation of energy/mass.
\end{rem}
\section{Applications}
\subsection{Linear estimates on $\R^n$}
In the specific case of the linear equation ($\eps=0$), one may extend
the identities of Theorems \ref{t1} and \ref{t2} through a limiting
argument in the spirit of \cite{LuisNicola}. Theorem \ref{bornebilineaire}
is in fact a consequence of upcoming Theorem \ref{t1l}. We start with the $1D$
case, which can also be derived by an explicit computation in Fourier
space, see \cite{OzTsu}. We will use the following definition of  the Fourier transform
of a function f 
$$\hat f(\xi)=\int_{\R^n} e^{- i2\pi x\xi} f(x)\,dx.$$
\begin{theoreme}[Ozawa-Tsutsumi \cite{OzTsu}]
  \label{t2l}
Let $n=1$, $u,v$ two solutions to \eqref{equ}, \eqref{eqv} with $\eps=0$, then
\begin{equation}
  \label{eq:lid1D}
\int_{\R\times\R}|\partial_x (u \bar v)|^2 \, dx dt =
4\pi\int_{\R\times\R} |\xi-\eta| |\hat u_0|^2(\xi)|\hat
v_0|^2(\eta)\,d\xi d\eta.
\end{equation}
\end{theoreme}
In higher dimensions, one has
\begin{theoreme}
  \label{t1l}
Let $\omega \in \R^n$ with $|\omega|=1$, $u$ solution to
\eqref{equd} with $\eps=0$. Then, with $x=x^\perp+s\omega$
\begin{multline}
\int_s |\partial_s (R(|u|^2))(s,\omega)|^2\,ds\\
{}+\int_t\int_s \int_{x\cdot \omega=s}\int_{y\cdot \omega=s} |u(x^\perp+s\omega)\partial_{s}
u(y^\perp+s \omega)-u(y^\perp+s\omega)\partial_{s}
u(x^\perp+s \omega)|^2 \,dx^\perp dy^\perp  dsdt\\
=4\pi\int_{\R^n\times\R^n} |\omega\cdot(\xi-\eta)| |\hat u_0|^2(\xi)|\hat
u_0|^2(\eta)\,d\xi d\eta.
\end{multline}
\end{theoreme}
Theorem \ref{t1} may be used in a different direction, recovering a
known bound for the linear equation (see \cite{Bourgain98}).
\begin{proposition}
  Let $u$ and $v$ be two solutions to \eqref{equ}, with
  $\eps=0$ and data $u_0,v_0$. Assume moreover that $\supp \hat u(\xi) \subset
  \{|\xi|\leq 2^k\}$ and $\supp \hat v(\xi-\xi_0) \subset
  \{|\xi|\leq 2^k\}$, with $|\xi_0|\sim 2^j$ and $k<<j$ (hence, the
  Fourier supports are separated and at distance roughly $2^j$). Then
  \begin{equation}
\label{eq:Bourgain}
\|u v\|^2_{L^2_{t,x}} \lesssim 2^{(n-1)k-j} \| u_0\|^2_{L^2_x}\|v_0\|^2_{L^2_x}.
  \end{equation}
\end{proposition}
\subsection{Scattering in $\R^n$}
A simple application of Theorems \ref{t1} and \ref{t2} is to recover
Nakanishi's scattering result for the $H^1$-subcritical (and
$L^2$-supercritical) defocusing equation. Such an alternative proof
was mentioned in \cite{CKSTTmora} for the cubic defocusing NLS in $3D$
(the authors actually proved a better result, as scattering is proved
to hold for $H^s$, $s>4/5$), and done in detail for the aforementioned
range in $1D$ in \cite{CHV}, where an a priori $L^8_{t,x}$ bound was
derived from a four particles interaction Morawetz inequality.
\begin{theoreme}\label{thscattering}
  Let $u_0\in H^1(\R^n)$, $n\geq1 $ and $u$ be the associated solution
  to \eqref{equ} with $\eps=1$, $1+\frac 4 n<p<1+
\frac{4}{n-2}$. Then one has scattering and
  polynomial bounds on space-time norms in term of mass $M=\int |u_0|^2$
  and energy $E=\int |\nabla_x u_0|^2+ 2 |u_0|^{p+1}/(p+1)$.
\end{theoreme}
 \begin{rem}
   Theorem \ref{thscattering} may also be found in \cite{CGTmora}. In
   fact, the authors go beyond the $H^1$ theory and establish
global existence and scattering for $H^s$ data, where $s_p<s<1$ and
$s_p$ is a critical exponent for their argument. One may also consult
the very recent survey \cite{GVsurvey}, which encompasses all the
known results, including ours, as well as extends the argument to
Hartree equations.
 \end{rem}

\subsection{Existence and scattering on a $3D$ exterior domain}
\label{sec:esext}
Due to the unavailability of scale-invariant
Strichartz estimates, the scattering issue is more difficult. In
fact, well-posedness in the energy class is already a significantly
more difficult problem, and is known to hold up to $p<3$
(\cite{BGTexterieur}), $p=3$ (\cite{Oana} and \cite{Ramonaext}); in
these references, non sharp (non scale-invariant) Strichartz estimates are obtained and
turned into the local existence result. Note that in
$3D$, our estimate is better with respect to scaling (sharp estimate with a
loss of a $1/4$ derivative) but somehow restrictive due to both the time
integrability range and the derivative loss. We first deal with
existence.
\begin{theoreme}
  \label{existence3D}
Let $1< p< 5$ and $n=3$. Let $u_0\in H^1_0(\Omega)$, with $\Omega$ an exterior
domain where local smoothing holds, $K$ a compact set such that
$\Omega^c\subset\subset K$. Then there exists a local in time solution $u$ to
\eqref{equd} which is $C_t(H^1_0(\Omega))$. Uniqueness holds in
$C_T(H^1_0)\cap L^4_T(W^{\frac 3 4,4})\cap  L^2_T(H^\frac 3 2(K))\cap L^4_T(L^\infty_x(K^c))$. Moreover, when $\eps=1$ (defocusing
case), the solution is global in time.
\end{theoreme}
When the domain $\Omega$ is star-shaped, one may use Proposition
\ref{p2d} and use the same strategy as in the $\R^n$ case to obtain scattering for the cubic
equation (with some significant additional technical difficulties, due
to the lack of the full set of Strichartz estimates). 
\begin{theoreme}
  \label{scattering3D}
Let $p=3$. Let $u_0\in H^1_0(\Omega)$, with $\Omega$ the exterior
 of a star-shaped domain. Then the global in time solution $u$ to
the defocusing equation \eqref{equd} scatters in $H^1_0(\Omega)$.
\end{theoreme}
\section{Proofs and further developments}

\subsection{1D computation, nonlinear equation}
\label{sec:1dle}
As a warm-up for subsequent computations, we prove Theorem \ref{t2} in
the special case $u=v$. Let $u$ be a solution to \eqref{equ}, and let
\begin{equation}
  I=\int_{x>y} (x-y) |u(x)|^2|u(y)|^2 \,dx\,dy.
\end{equation}
Compute the time derivative of $I$: as we have $i\partial_t u+\partial_x^2 u=\eps |u|^{p-1} u=f$, the
nonlinear part vanishes when computing
\begin{equation}
 \frac{ d |u(x)|^2}{dt}=\frac 1 i ( u \partial^2_x \bar u-\bar
 u \partial_x^2 u)= i \partial_x (\bar u\partial_x u-u\partial_x \bar
 u)=-2\partial_x(\Im(\bar u \partial_x u)),
\end{equation}
and we have
\begin{eqnarray*}
  \partial_t I & = & -2 \int_{x>y} (x-y) \left(  \partial_x (\Im(\bar u \partial_x u))(x)|u(y)|^2+  \partial_y ( \Im(\bar u \partial_y u))(y)|u(x)|^2 \right)\,dx\,dy\\
&  = & 2  \left (\int_{x>y}   \Im(\bar u
  \partial_x u)(x)|u(y)|^2- \Im( \bar u
  \partial_y u)(y)|u(x)|^2 \,dy\,dx\right).
\end{eqnarray*}
Derive again in time and focus on the first term: it will be a sum of
3 terms ($K_1,K_2,K_3$),
\begin{eqnarray*}
    K_1 &  = &  2 \int_{x>y}   \Im(\bar u
  \partial_x u)(x)  \frac{d |u(y)|^2}{dt} \,dy\,dx =  -4 \int_{x>y}   \Im(\bar u
  \partial_x u)(x)  \partial_y \Im ( \bar u \partial_y u)(y)\,dy\,dx\\
   & = &-4 \int_x (\Im (\bar u\partial_x u))^2(x) \,dx.
\end{eqnarray*}
\begin{rem}
  Notice for further use that when picking the second term in $\dot
  I$, it will contribute exactly another $K_1$ term (boundary term
  with opposite sign).
\end{rem}
Now, the second term is the sum of a linear term,
\begin{eqnarray*}
  K_2  & = &  \int_{x>y}   ( i \partial_t u \partial_x \bar u- i \partial_t \bar u
  \partial_x u)(x)|u(y)|^2 \,dy\,dx   =  \int_{x>y}   ( - \partial_x (|\partial_x u|^2)(x)|u(y)|^2
  \,dy\,dx\\   &  = & \int_{y}   |\partial_y u|^2(y)|u(y)|^2 \,dy =
  \int_{y}   |\bar u \partial_y u|^2(y)\,dy,
\end{eqnarray*}
and a nonlinear term,
\begin{eqnarray*}
  A_2  &   = & \int_{x>y}   ( f \partial_x \bar u+ \bar f
  \partial_x u)(x)|u(y)|^2 \,dy\,dx = \int_{x>y}    |u|^{p-1}(x) \partial_x (|u|^2)(x)|u(y)|^2
  \,dy\,dx\\
 &   = & -\frac 2 {p+1}\int_{y} |u|^{p+3}(y) \,dy.
\end{eqnarray*}
The same remark applies for the other contribution with $x$ and $y$
reversed (so we double $K_2+A_2$). The next term is,
\begin{eqnarray*}
  K_3  & = &  \int_{x>y}   ( i u \partial_x\partial_t \bar u- i \bar u
  \partial_x\partial_t u)(x)|u(y)|^2 \,dy\,dx\\
 &   =  &\int_{x>y}   -( i \partial_x u \partial_t \bar u- i\partial_x \bar u
  \partial_t u)(x)|u(y)|^2 \,dy\,dx+K_4 = K_2+K_4,
\end{eqnarray*}
with $K_4$ being the boundary term, namely
\begin{eqnarray*}
    K_4 &  = & - \int_{y}   ( i u \partial_t \bar u- i \bar u
  \partial_t u)(y)|u(y)|^2 \,dy  =  - \int_{y}   (  u \partial^2_y \bar u+ \bar u
  \partial^2_y u)(y)|u(y)|^2 \,dy\\
     &  = &   \int_{y}   2 |\partial_y  u|^2(y)|u(y)|^2+(u\partial_y \bar u+\bar u\partial_y u)\partial_y
  (|u|^2) \,dy = 2 K_2 + \int_y (\partial_y
  (|u|^2))^2 \,dy.
\end{eqnarray*}
So that
\begin{equation*}
K_2+K_3 =  4  \int_{y}   |\partial_y u|^2)(y)|u(y)|^2 \,dy+ \int_y (\partial_y
  (|u|^2))^2 \,dy.
\end{equation*}
The nonlinear contribution $A_3$ verifies the same identity, namely
$A_3=A_2+A_4$, and $A_4$ is the following nonlinear boundary term:
\begin{equation*}
      A_4   =   \int_{y}   (  u \bar f+ \bar u f )(y)|u(y)|^2 \,dy =   2 \int_{y} |u(y)|^{p+3} \,dy
\end{equation*}
and the total contribution of the nonlinear term is
\begin{equation*}
  A = A_2+A_3=  (2- \frac{4}{p+1})\int_y |u|^{p+3}(y) \,dy.
\end{equation*}
Now the claim is that the second part in $\partial_t I$ gives the exact same expression: $x$ and $y$ are exchanged, we have a
  minus sign in front, and the boundary term will be at the opposite
  end, switching the sign. Hence,
\begin{eqnarray*}
  \frac{d^2 I}{d t^2} &  = &  2(K_1+K_2+K_3+A) \\
  \frac{d^2 I}{d t^2}-2A &   = & 2 \left (4 \int_x (\Re (u\partial_x \bar u))^2(x) \,dx+  \int_y (\partial_y
  (|u|^2))^2 \,dy\right)\\
  \frac{d^2 I}{d t^2} &   = & 4 \int_x (\partial_x (|u|^2))^2(x) \,dx+ 4 \eps \int_y |u|^{p+3}(y) (1- \frac{2}{p+1}) \,dy.
\end{eqnarray*}
which is nothing but the identity in Theorem \ref{t2}. Notice that $I$
is a convex function whenever $\eps=0,1$.
\begin{rem}
  One may somewhat shorten the proof by introducing the density of
  mass $N$, the current $J$ and the (one dimensional for now !) ``tensor'' $T$,
$$
N=|u|^2,\,J=2\Im(\bar u \partial_x u),\,T=4 |\partial_x u|^2-\Delta
N+\varepsilon(2-\frac 4 {p+1}) N^{\frac{p+1}2},
$$
and then use local conservation laws to perform the integrations by parts,
$$
\partial_t N+\partial_x J=0 \,\,\text{   and   }\,\,\partial_t J+\partial_x T=0.
$$
Evidently, the relation $\partial_t^2 N=\partial^2_x T$ is
behind any sort of virial identity, bilinear or not, and the reader
may consult \cite{GVsurvey} for a very nice survey of bilinear virial
estimates, including ours, which present the above derivation in a
concise and elegant form.
\end{rem}
\subsection{A digression on local smoothing estimates}
\label{localsmoothingsection}
The 1D proof from the previous section makes crucial use of boundary
terms $x=y$ arising in integrations by parts. Let us now give an
elementary proof of the following well-known 1D estimate
(\cite{KPVlequel}).
\begin{proposition}
  Let $u$ be a solution to the linear Schr\"odinger equation on $\R$:
  \begin{equation}
    \label{eq:smoothing1D}
    \sup_x \int_\R |\partial_x u|^2 (x,t)\,dt = C \|u_0\|^2_{\dot
    H^\frac 1 2}.
  \end{equation}
\end{proposition}
Consider $v(x)=u(x)-u(-x)$ the odd part of $u$: $v$ still satisfies
the Schr\"odinger equation (in fact, $v$ may be seen as a solution to
the equation on $\R_+$ with Dirichlet boundary condition
$v(x=0)=0$). Multiply the equation for $v$ by $\partial_x \bar v$ and
integrate between $x=\infty$ and $x=0$:
\begin{equation*}
  \int_{t_1}^{t_2} \int_\infty^0 i(\partial_t v \partial_x \bar v -\partial_t \bar v
  \partial_x v)+ \int_{t_1}^{t_2} |\partial_x v|^2(0)  =  0.
\end{equation*}
A double integration by parts in the first term yields two different
types of boundary terms: time slice ones,
$$
|\int_\infty^0 \Im \bar v \partial_x v \,dx (t_1)-\int_\infty^0 \Im
\bar v \partial_x v \,dx (t_2)|\lesssim \sup_{[t_1,t_2]}\|v\|^2_{\dot
  H^\frac 1 2},
$$
where one is using duality and $v(x=0)=0$. On the other hand, one has
a remaining spatial boundary term,
$$
\int_{t_1}^{t_2} \Im \bar v \partial_t v (0)=0
$$ 
due to the boundary condition. Inequality in \eqref{eq:smoothing1D}
follows trivially by translation invariance, noticing that $\partial_x
v(0)=2\partial_x u(0)$. Sending both $t_1$ and $t_2$ to $\pm\infty$ and
recalling the asymptotic of the free solution would provide the
equality by rewriting the momentum in term of $\hat u_0$ (see \cite{LuisNicola}).

Alternatively, one may derive this estimate by computing twice the
time derivative of
$$
I=\int_{x>y} (x-y) |v_y|^2(x,t)\,dx, \text{ with } v_y=u(x+y)-u(y-x).
$$
If one goes to dimension $n$, we may instead consider the reflexion
with respect to the hyperplane $x_n=0$, compute
$$
I=\int_{x_n >y_n } (x_n-y_n) |v_y|^2(x,t)\,dx, \text{ with } v_y(x)=u(x',x_n+y_n)-u(x',y_n-x_n).
$$
In the computation, one may pick up additional boundary terms, namely
$\int_t |\nabla' v_y|^2(x',y_n)\,dx' dt$, which vanish thanks to
$v_y(x_n=y_n)=0$. Hence, we have obtained a very elementary proof of
the following variant of the local smoothing effect, with no use of
the Fourier transform in space or time. One, however, relies heavily
on the invariances.
\begin{proposition}
  Let $u$ be a solution to the linear Schr\"odinger equation on
  $\R^n$, and $\omega$ a direction, with $x=(x_\omega^\perp,x_\omega)$:
  \begin{equation}
    \label{eq:smoothingnD}
    \sup_{x_\omega} \int_{\R\times \R^{n-1}} |\partial_{x_\omega} u|^2
    (x,t)\,dx^{\perp}_\omega dt \lesssim \|u_0\|^2_{\dot
    H^\frac 1 2}.
  \end{equation}
\end{proposition}
In view of this computation, the weight $\rho(x-y)=(x-y)\cdot \omega$
appears to be a rather natural choice in $\R^n$, when trying to
average the virial on the half-space $(x-y)\cdot \omega>0$.
\subsection{Bilinear estimate on the nonlinear equation, the general
  case}
We now turn our attention to the general case, and prove Theorems
\ref{t1},\ref{t2}, \ref{t3}, \ref{t1d} all together. We consider the equation on a domain $\Omega$, with Dirichlet boundary conditions $u_{|\partial \Omega =0}$. Recall that
\begin{equation}
\label{dtM}
    i\partial_t (|u|^2)  =  u\Delta \bar u-\bar u \Delta u  =  \nabla
    \cdot (u\nabla \bar u -\bar u \nabla u)=-2i \nabla\cdot \Im (\bar u \nabla u).
\end{equation}
Set
\begin{equation}
  I =\int_{\Omega\times \Omega} \rho(x-y) |u|^2(x)|v|^2(y) \,dxdy.
\end{equation}
We compute
\begin{eqnarray*}
  \partial_t I & = & - 2 \int \rho \left(|v|^2 \nabla \cdot \Im (\bar u
    \nabla u)+|u|^2 \nabla \cdot \Im (\bar v \nabla
  v)\right)\\
    & = &   2 \int \nabla_x \rho \cdot \left(|v|^2(y) \Im (\bar u \nabla
  u) (x)-|u|^2(x) \Im (\bar v \nabla
  v )(y)\right) \, dxdy.
\end{eqnarray*}
where there is no boundary term when applying Stokes, as there is
always a factor of $u$ or $v$ to cancel such a term due to the
Dirichlet condition, and we used $\nabla_x \rho =- \nabla_y \rho$.

Now, we compute $\partial^2_t I=J_x+J_y+J_{xy}$ depending on where the
time derivative lands (with obvious notations). We have
\begin{equation}
  J_x = \int |v|^2(y) \nabla_x \rho \cdot \partial_t  \left(\frac{\bar u \nabla
  u -u \nabla \bar u}{i}\right)(x) \,dxdy.
\end{equation}
Now
\begin{eqnarray*}
\partial_t   \left(\frac{\bar u \nabla
  u -u \nabla \bar u}{i}\right) & = & (-\Delta u+\eps |u|^{p-1}
  u)\nabla \bar u +(-\Delta \bar u+\eps |u|^{p-1}
  \bar u)\nabla  u  \\
 & & {}- \left[ u \nabla(-\Delta \bar u+\eps |u|^{p-1}
  \bar u)+ \bar u \nabla(-\Delta  u+\eps |u|^{p-1}
   u)\right]\\
 & = & -\Delta u \nabla \bar u-\Delta \bar u \nabla u+u \nabla \Delta \bar u+\bar u \nabla \Delta u-\eps |u|^2 \nabla (|u|^{p-1}).
\end{eqnarray*}
Back to $J_x$, we call $K_1, K_2$  the bilinear and nonlinear terms coming from the
above formula. We use Einstein convention for summation :
\begin{eqnarray*}
  K_1 &  = & \int |v|^2 \partial_i \rho (-\partial_k \partial^k u \partial^i
  \bar u -\partial_k \partial^k \bar u \partial^i u)+\int
  |v|^2 \partial_i \rho(u \partial^i
  \partial_k \partial^k \bar u +\bar u \partial^i
  \partial_k \partial^k u)\\
 & = & K_{11}+K_{12}.
\end{eqnarray*}
We have ($n(x)$ being the outgoing normal vector at $x\in \partial\Omega$)
\begin{eqnarray*}
  K_{11} & = & \int |v|^2(y) \nabla_x \rho(x-y) \cdot (\nabla_x \bar u (-\nabla_x\cdot \nabla_x u )+ \nabla_x u(-\nabla_x\cdot \nabla_x \bar u))(x)\\
 & = & - \int |v|^2(y) (\nabla_x \rho\cdot \nabla_x \bar u \nabla_x u\cdot n(x) +\nabla_x \rho\cdot \nabla_x  u \nabla_x \bar u\cdot n(x)) \,dS_x\\
& & {}+\int |v|^2 (y) (\nabla_x u\cdot \nabla_x (\nabla_x \rho\cdot \nabla_x \bar u+\nabla_x \bar u\cdot \nabla_x (\nabla_x \rho\cdot \nabla_x  u)(x) \\
 & = & - 2 \int |v|^2(y) \partial_n \rho(x-y) |\partial_n  u|^2 (x) \,dS_x \,dy\\
 &  & +{} 2 \int |v|^2 \partial_i \partial^k \rho\, \partial_k u \partial^i
  \bar u \text{     (recall the Hessian is symmetric)} \\
 & & {}+ \int |v|^2 \partial_i \rho\,(\partial_k u \partial^i \partial^k
  \bar u+\partial_k \bar u \partial^k \partial^i u)
\end{eqnarray*}
where we used the Dirichlet condition in the boundary term and
expanded the remaining terms. On the other hand, as all boundary terms cancel due to the Dirichlet condition,
$$
 K_{12}  = -\int |v|^2 \partial_i \partial^k \rho\, (u\partial^i \partial_k
  \bar u +\bar u\partial^i \partial_k u) 
 - \int |v|^2 \partial_i \rho \,(\partial^k u \partial^i
  \partial_k \bar u+\partial^k \bar u \partial^i \partial_k u).
$$
Summing $K_{11}$ and $K_{12}$, their respective last terms  cancel each other. Integrate the first term in $K_{12}$ with respect to $\partial^i$, there is (again) no boundary term, and finally
\begin{multline*}
      K_1   =  4 \int |v|^2(y) \partial_i \partial^k \rho(x-y)\, \partial_k u \partial^i
  \bar u(x) +\int |v|^2(y) \partial^i \partial_i \partial^k
  \rho(x-y)\, \partial_k( |u|^2)(x) \, dxdy \\
{}- 2 \int |v|^2(y) \partial_n \rho(x-y) |\partial_n  u|^2 (x) \,dS_x \,dy
\end{multline*}
However, $\partial^2_{x_i} \partial_{x_k}
\rho=-\partial^2_{x_i}\partial_{y_k} \rho=-\partial_{x_i}\partial_{x_k}\partial_{y_i}$, so that one
last integration by parts yields, denoting by $H_\rho$ the Hessian of $\rho$,
\begin{multline*}
    K_1= 4 \int |v|^2(y) H_\rho(x-y)(\nabla_x u(x),\nabla_x \bar u(x))
  \,dxdy+ \int H_\rho(x-y) (\nabla_x(|u|^2)(x), \nabla_y(|v|^2)(y))\,dxdy\\
{}- 2 \int |v|^2(y) \partial_n \rho(x-y) |\partial_n  u|^2 (x) \,dS_x \,dy,
\end{multline*}
given that the integration by parts in $y$ does not have a boundary
term either, and where me may freely replace the second term using the identity
\begin{equation*}
   \int H_\rho(x-y) (\nabla_x(|u|^2)(x), \nabla_y(|v|^2)(y))\,dxdy=
   \int \Delta\rho(x-y) \nabla_x(|u|^2)(x)\cdot  \nabla_y(|v|^2)(y)\,dxdy.
\end{equation*}

Now, we go back to the nonlinear term $K_2$:
\begin{eqnarray*}
  K_2 & = & -\int |v|^2(y) \nabla_x \rho \cdot \eps (|u|^{p-1})^{\frac
  2 {p-1}} \nabla (|u|^{p-1}) \,dxdy\\
   & = & - \eps \int |v|^2(y) \nabla_x \rho \cdot \nabla
  (|u|^{p+1})\frac{1}{\frac 2 {p-1}+1} \,dxdy\\
 & = & \eps \int |v|^2(y) (\Delta_x \rho)(x-y)  |u|^{p+1} (x) \frac{p-1}{p+1} \,dxdy,
\end{eqnarray*}
performing one more integration by parts (with no boundary  term !). Assuming that $\rho(x-y)=\rho(y-x)$, the second term $J_y$ is exactly $J_x$ by
symmetry, up to permutation of $u$ and $v$. We are left with
\begin{equation*}
  J_{xy}= J_1+J_2,
\end{equation*}
where again by symmetry both terms are equal and
\begin{eqnarray*}
  J_1 & = & 2\int \Im(\bar u \nabla u)(x) \cdot \nabla_x \rho
  \partial_t (|v|^2(y))\,dxdy\\
 & = & 4 \int \Im(\bar u \nabla u)(x) \cdot \nabla_x \rho
  \nabla_y \Im(\bar v \nabla v(y))\,dxdy
\end{eqnarray*}
and using again $\nabla_y \rho=-\nabla_x \rho$, we integrate by parts in $y$ (with no boundary term)
$$
  J_1  =  -4  \int H_\rho(x-y) \left(\Im(\bar u \nabla u)(x),\Im(\bar v
  \nabla v)(y)\right) \,dxdy.
$$
Finally,
\begin{align*}
  \partial^2_t I  = & 4 \int \left( |v|^2(y) H_\rho(x-y)(\nabla u(x),\nabla \bar u(x))
   + |u|^2(x) H_\rho(x-y)(\nabla v(y),\nabla \bar v(y))
  \right)\,dxdy\\
  & {}+2\int \left( H_\rho(x-y) (\nabla(|u|^2)(x),\nabla(|v|^2)(y))-4 H_\rho(x-y) \left(\Im(\bar u \nabla u)(x),\Im(\bar v \nabla v)(y)\right)\right) \,dxdy\\
  & {}+ \eps (1-\frac{2}{p+1})\int \left(|v|^2(y) (\Delta\rho)(x-y)  |u|^{p+1} (x)
  + |u|^2(x) (\Delta_x \rho)(x-y)  |v|^{p+1} (y)\right) \,dxdy\\
 & {}-2 \int |v|^2(y) \partial_n \rho(x-y) |\partial_n  u|^2 (x) \,dS_x \,dy-2 \int |u|^2(x) \partial_n \rho(x-y) |\partial_n  v|^2 (y) \,dS_y \,dx.
\end{align*}
By definition, $H_\rho$ is symmetric. Then one may diagonalize and be
left with just one direction (or, more accurately, a diagonalized
matrix). Discarding a factor 2 and the eigenvalue $\lambda(x-y)$, we set (where $\partial$ denotes derivation in the direction of
the eigenvector associated to $\lambda$)
\begin{multline*}
 \Gamma  =  2|v|^2(y) |\partial u|^2(x)+ 2|u|^2(x) |\partial v|^2(y)\\{}+ (v\partial\bar v+\bar v \partial v)(y) (u\partial\bar
       u+\bar u \partial u)(x)- (v\partial\bar v-\bar v \partial v)(y) (\bar
       u\partial u-  u \partial\bar u)(x).
\end{multline*}
Expanding the last two terms and canceling out, we get
\begin{eqnarray*}
 \Gamma & = & 2|v|^2(y) |\partial u|^2(x)+ 2|u|^2(x) |\partial
 v|^2(y)+ 2 v\partial\bar v(y) u\partial\bar u(x) +2 \bar v \partial
 v(y)\bar u\partial u(x)\\
& = &  2 |\bar v (y)\partial u(x)+u(x)\partial \bar v(y)|^2.
\end{eqnarray*}
Now, one may rewrite $\Gamma$ in a different way, by taking advantage of
the identity
$$
 |v (y) \partial u(x)-u(x)\partial v(y)|^2 + \partial (|v|^2)(y) \partial(|u|^2)(x)= |\bar v(y) \partial
   u(x)+u(x) \partial \bar v(y)|^2.
$$
This achieves the proof of the first part of Theorem \ref{t3}, namely
\eqref{nonlinearmauvais} and \eqref{nonlinearrho}. Finally, in the special case $n=1$ we
obtain Theorem \ref{t2}.  A straightforward generalization of the
$n=1$ case will follow by setting $\rho=|(x-y)\cdot \omega|$ with
$\omega\in \S^n$: we will obtain Theorem
\ref{t1} and Theorem \ref{t1d}.
\begin{rem}
  Notice that if $u=v$ and $x=y$, one recovers the same identity, with
  $\Gamma=2(\partial(|u|^2))^2$, with both expressions, which is
  consistent with our previous $1D$ computation.
\end{rem}
Set $\rho=|x_n-y_n|$ for convenience, and let us focus on the linear
equation in $\R^n$; we have obtained, discarding a positive term, the inequality
\begin{equation}
\label{discard}
\int_{-T}^T   \int_{x_n} \left(\partial_n\left(\int_{x'} |u|^2(x',x_n,t)\,dx'\right) \right)^2 \,dx_n dt \lesssim
\int_{x_n<y_n} \Im \bar u \partial_n u(x) |u|^2(y)\,dxdy |_{-T}^T.
\end{equation}
Proceeding exactly as in \cite{LuisNicola}, one may send $T\rightarrow
+\infty$ and recover an exact formula for the right handside. Recall the following asymptotic formula
for the solution $U(t, z)$ to the linear equation $i\partial_t
U+\Delta U=0$ with data $U_0$, and $z\in \R^{m}$, which follows
directly from the explicit representation as a convolution by the Gaussian
kernel $(4\pi i t)^{-m/2}\exp(i|z|^2/4t)$:
\begin{equation}\label{asymptotique+}
\lim_{t\rightarrow \pm \infty}  \|U(t,z) - 
\frac{e^{\pm {i} \frac{|z|^2}
{4t}}}{(4\pi i t)^{m/2}} 
\hat U_0 \left( \pm \frac z {4\pi|t|} \right)
\|_{L^2({\mathbb R^m})}=0.
\end{equation}
By using \eqref{asymptotique+} we get
$$\lim_{t\rightarrow +\infty} \|U(t,z)-V(t,z)
\|_{L^2({\mathbb R}^m)}=0,$$
where $V(t,z):=e^{-{i}m\pi/4}\frac{e^{{i} \frac{|z|^2}{4t}}}
{(4\pi t)^{m/2}} 
\hat U_0 \left( \frac z{ t} \right)$. On the other hand, for any
direction $s$, $\partial_s
U$ is also a solution, hence
$$\lim_{t\rightarrow +\infty} \|\partial_s U(t, z)- 
W_s(t, Z)\|_{L^2(\mathbb R^m)}=0$$
where
$$W_s (t, z):= e^{-{i}m\pi/4}\frac{e^{{i} 
\frac{|z|^2}{4t}}}{(4\pi t)^{m/2}} 
{i} \frac{s}{2 t}\hat U_0\left( \frac z{4\pi t} \right).$$
We easily deduce 
\begin{equation*}\lim_{t\rightarrow +\infty}
\int_{{\mathbb R}^m} [\bar U(t, z) \partial_s  U(t, z) 
- \bar V(t, z)  W_s(t,z)] \phi(s)dz=0
\end{equation*}
for any $\phi\in L^\infty$.  Then if $\phi(s)=\partial_s |s|$,
\begin{align*}
  \lim_{t\rightarrow +\infty} {\Im}\int
_{{\mathbb R}^m} \bar U(t, z) \phi(s) \partial_s U(t, z) dz
 & =\lim_{t\rightarrow +\infty}
{\Im} \int_{{\mathbb R}^m} \bar V(t,z) 
W_s(t, z) \phi(s) dz  \\
 & =\lim_{t\rightarrow +\infty}
{(4\pi t)^{-m} }
\int_{{\mathbb R}^m} \frac{s}{2t} \left
|\hat U_0 \left ( \frac z{ 4\pi t}\right )
\right |^2\phi( s)dz\\
 & =4\pi \int_{{\mathbb R}^m} \frac{|s|}{2} \left
|\hat U_0 ( z)\right |^2dz.
\end{align*}
Let $u$ be a solution of the linear equation in $\R^n$. We proceed
with a tensor product solution and set $U(t,z)=u(t,x)u(t,y)$ with $z=(x,y)$, $m=2n$, and pick the
direction $s=x_n-y_n$. Then the limit when $T\rightarrow +\infty$ of
the right handside in \eqref{discard} will be a multiple of $\int
|x_n-y_n| |\hat u_0|^2(x)|\hat u_0|^2(y)\,dxdy$. Hence we have obtained
\begin{equation}
\int_{t}   \int_{x_n} \left(\partial_n\left(\int_{x'} |u|^2(x',x_n,t)\,dx'\right) \right)^2 \,dx_n dt \lesssim
\int|\xi_n-\eta_n| |\hat u_0|^2(\xi) |\hat u_0|^2(\eta)\,d\xi d\eta.
\end{equation}
Up to the upcoming introduction of the Radon transform, this is exactly Theorem \ref{bornebilineaire} but with $u=v$. Applying
this estimate to $u+v$ and $u+iv$ allows to control both $|\partial_n
(\Re(u\bar v))|$ and $|\partial_n (\Im(u\bar v))|$, and applying
Cauchy-Schwarz repeatedly on the right-hand side, we obtain Theorem
\ref{bornebilineaire}. On the other hand, in the special case $n=1$,
one does not need to set $u=v$ and we obtain Theorem
\ref{t2l}. Finally, if we retain the discarded term in \eqref{discard}
and keep $u=v$, we obtain Theorem \ref{t1l}.

Dilating $u$ and $v$ in opposite way and optimizing allows us to
replace the right-hand side by $\|u_0\|_{L^2_x} \|v_0\|_{\dot H^\frac 1
  2}+ \|v_0\|_{L^2_x} \|u_0\|_{\dot H^\frac 1 2}$, up to Proposition
\ref{propH12} whose proof we postpone for the moment.

Now we reduce the directional estimate we obtained to a generic one by
introducing the Radon transform. By rotation, one may replace $x_n$ by the coordinate along any
direction $\omega$, so that if $R(f)(s,\omega)$ is the Radon transform
of a function $f$, namely
\begin{equation}
  R(f)(s,\omega)=\int_{x\cdot \omega=s} f \,d\mu_{s,\omega},
\end{equation}
where $\mu_{s,\omega}$ is the induced measure on the hyperplane
$x\cdot\omega=s$, the previous estimate can be recast as
\begin{equation}
\sup_\omega \int_t \int_s |\partial_s (R(u\bar v))(s,\omega)|^2\,ds dt \lesssim \|u_0\|^2_{L^2_x} \|v_0\|_{\dot H^\frac 1
  2}^2+ \|v_0\|^2_{L^2_x} \|u_0\|_{\dot H^\frac 1 2}^2.
\end{equation}
Replacing the $L^\infty_\omega$ by $L^2_\omega$ and using that
\begin{equation}
\label{radonplancherel}
  \| |\partial_s|^{\frac{n-1}2} R(f)\|_{L^2}=\|f\|_{L^2},
\end{equation}
one recovers known bounds on the linear equation for $n=2$ and $n=3$:
\begin{itemize}
\item if $n=2$,
  \begin{equation}
\label{L4H12}
   \int_{t_1}^{t_2} \||\nabla|^{\frac 1 2} (u\bar
   v)\|^2_{L^2}\,dt\lesssim \|u_0\|^2_{L^2_x} \|v_0\|^2_{\dot H^\frac 1
  2}+ \|v_0\|^2_{L^2_x} \|u_0\|_{\dot H^\frac 1 2}^2.
  \end{equation}
One may get an $L^4_t(L^8_x)$ bound for $\dot H^{\frac 1 4}$ data, but fails short of getting the usual $L^4_{t,x}$ bound. We will recover this bound through a refined analysis using the Radon bound in a more efficient way.
\item if $n=3$
  \begin{equation}
   \int_{t_1}^{t_2} \|(u\bar
   v)\|^2_{L^2}\,dt\lesssim \|u_0\|^2_{L^2_x} \|v_0\|^2_{\dot H^\frac 1
  2}+ \|v_0\|^2_{L^2_x} \|u_0\|_{\dot H^\frac 1 2}^2
  \end{equation}
which is a (linear) variation on the original $L^4_{t,x}$ Morawetz
interaction estimate from \cite{CKSTTmora}.
\end{itemize}
Next, one would like to take advantage of the $L^\infty_\omega$
bound. Consider the situation where $v$ is frequency localized in a
(small) ball $|\xi|\lesssim 2^k$ and $u$ is frequency localized in a
ball of the same size but which is included in the annulus $|\xi|\sim
2^j$, with $k<<j$. The Fourier transform of $u\bar v$ has roughly the
same frequency localization as $u$, hence it is supported in a ball of
size $2^k$ and in an angular sector of (angular) width $2^{(n-1)(k-j)}$
(the volume of the $(n-1)$ dimensional cap which is the intersection of the angular sector and the sphere of radius $1$). As the
Fourier transform of $R(f)$ is connected with $\hat f$ by the
following formula,
\begin{equation}\label{radonfourier}
g(\rho,\omega)={\mathcal F}_{s\rightarrow \rho} (Rf(., \omega))[\rho]=
\hat f(\rho\omega) \hbox{ } \forall \omega \in                                            
{\mathbb S}^{n-1},
\end{equation}
we will have, for such $f=u\bar v$, by Plancherel,
\begin{equation*}
\int_{\rho,\omega} \rho^{3-n} |\hat f(\rho\omega)|^2
\rho^{n-1}\,
d\rho d\omega\lesssim
  2^{(n-1)(k-j)} \sup_\omega \int_\rho |\rho g(\rho,\omega)|^2 \,d\rho
\end{equation*}
which translates into
\begin{equation*}
  \||\nabla|^\frac 1 2 (u \bar v)\|^2_{L^2_{t,x}}\lesssim 2^{(n-1)k-j} ( \|u_0\|_{\dot H^\frac 1 2}^2 \|v_0\|^2_{L^2}+\|v_0\|_{\dot H^\frac 1 2}^2 \|u_0\|_{L^2}^2),
\end{equation*}
and due to the frequency localization, one may remove the
half-derivative on both sides; the complex conjugate is now
irrelevant, and we get \eqref{eq:Bourgain}. For example
Bourgain's original estimate for $n=2$ (\cite{Bourgain98})reads
\begin{equation}
  \|u v\|^2_{L^2_{t,x}}\lesssim 2^{k-j}  \|u_0\|^2_{L^2}\|v_0\|_{L^2}^2.
\end{equation}
By a Galilean transform, one may shift both factors by any
$\xi_0$ in frequency space, as both norms on the right and the left
are galilean invariant. Thus, we obtain that for $u,v$ such that their
Fourier supports are in balls of size $2^{2k}$ which are $2^j$ apart,
\begin{equation}
\label{Patrick}
\|uv\|^2_{L^2_{t,x}} \lesssim 2^{k-j}  \|u_0\|^2_{L^2}\|v_0\|_{L^2}^2.
\end{equation}
Assuming only \eqref{Patrick}, one may then recover the usual $L^4_{t,x}$ bound by the usual Whitney
decomposition trick, see \cite{TVV}. However, we may derive it
directly: consider 
$$
|u|^2=\sum_j S_{j-2} \bar u \Delta_j u+\sum_j
S_{j-2} \bar u \Delta_j u+\sum_{|j-j'|\leq 1} \Delta_{j'} \bar u
  \Delta_j u,
$$
the usual paraproduct decomposition. On both paraproduct terms, we
take advantage of the frequency separation; applying \eqref{Patrick}
provides the $L^4_{t,x}$ bound. On the reminder term, we have to
consider (abusing notations by reducing the sum to the diagonal one)
$$
\Delta_k (\sum_{k\lesssim j} \Delta_j u \Delta_j \bar u).
$$
If $j>>k$, then only opposite balls of radius $2^k$ (and at distance
$2^j$ from $\xi =0$) contribute, and again we may use \eqref{Patrick}
and sum in $j$. When $k\sim j$, either the two supports are separated
and \eqref{Patrick} will do, or the supports are the same (splitting
in a finite number of smaller balls if necessary), but then
they do not overlap the origin in $\xi$: one may go back to
\eqref{L4H12} and take advantage of the support condition to get rid
of half a derivative.

Let us go back to the nonlinear equation: our choice of
$\rho(x-y)=|\omega\cdot (x-y)|$ in Theorem \ref{t3}, together with the
definition of the Radon transform, immediately yields Theorems
\ref{t1} and \ref{t1d}, as the former is a particular case of the
later. 

 We now prove Proposition \ref{propH12}, starting with \eqref{eq:moment}: but we almost did in section
\ref{localsmoothingsection}. At fixed $y_n$,
$$
\left |\int_{x_n<y_n} \Im  (\bar u(x)-\bar u(x',y_n)) \partial_n u(x) \,dx\right|
\lesssim \| u\|^2_{\dot H^\frac 1 2},
$$
using duality in $\dot H_x^{\frac 1 2}(\R^n_+)$, as $u(x)-u(x',y_n)\in \dot
H_x^{\frac 1 2}(\R^n_+)$. Then one may rewrite, 

\begin{multline*}
2 i \int_{x_n<y_n} |u|^2(y) \Im \bar u(x) \partial_n u(x)\,dx
dy = \int_{x_n<y_n} |u|^2(y) ( (\bar u(x)-\bar u(x',y_n)) \partial_n u(x)\\
{}-(u(x)-u(x',y_n)) \partial_n \bar u(x)) \,dx
dy,
\end{multline*}
as
$$
\int_{x_n<y_n}\partial_n u(x)\,dx_n=u(x',y_n) \imp \int_{x_n<y_n}\bar
u(x',y_n) \partial_n u(x)\,dx=\int_{x'} |u|^2(x',y_n) \,dx'.
$$
From there \eqref{eq:moment} easily follows.

We proceed with \eqref{eq:borneapriori}, which follows from averaging
 the Radon transform over directions $\omega$ in $L^2_\omega$ in \eqref{radonplancherel}: the linear part we already
obtained; now both $R(|u|^2)$ and
$R(|u|^{p+1})$ are positive, we immediately have by Cauchy-Schwarz
$$
|R(|u|^{\frac {p+3} 2})|^2\lesssim R(|u|^2)R(|u|^{p+1}).
$$
Discarding a positive term in the left handside of \eqref{idt1}, we finally obtain
\eqref{eq:borneapriori} which ends the proof.

\section{Local smoothing and control of the trace for NLS on a domain}
We now prove Proposition \ref{p1d}. Let us stress, once again, that
for the linear equation, \eqref{eq:nlsdnormal} holds on any
non-trapping domain for any dimension (see \cite{BurqDuke}). Hence, the
  purpose of this section is to provide a simple integration by parts
  proof when $n\geq 3$, which equally applies to the nonlinear
  defocusing equation. Let us consider again
$$
i\partial_t u+\Delta u-\eps |u|^{p-1} u=0,
$$
where $\Delta$ is the Laplacian with Dirichlet boundary condition
$u_{|\partial \Omega} =0$, and $\Omega$ is the exterior of a
star-shaped body with smooth boundary.

First, the virial identity (the following computation is standard and
we provide it for completeness): let us denote
\begin{equation}
  \label{eq:virielorig}
  M_h(t)= \int_\Omega |u|^2(x,t) h(x)\,dx,
\end{equation}
where $h$ is any smooth real-valued function on $\Omega$. Then compute
(recalling \eqref{dtM})
\begin{equation*}
  \frac{d}{dt} M_h(t) = -2 \Im\int h \nabla\cdot (\bar u \nabla u) = 2 \Im \int \bar u \nabla u \cdot \nabla h,
\end{equation*}
where we used the Dirichlet boundary condition when integrating by
parts. Now,
\begin{align*}
  \frac{d^2}{dt^2} M_h(t) & =  2 \Im \int (\partial_t \bar u \nabla
  u+\bar u \nabla \partial_t u) \cdot \nabla h = -2 \Im \int \partial_t u \left(2 \nabla
  \bar u \cdot \nabla h+\bar u\Delta h\right)\\
                          & = -2 \Re \int (\Delta u-\eps |u|^{p-1} u) \left(2 \nabla
  \bar u \cdot \nabla h+\bar u\Delta h\right)\\
                          & = -4 \Re \int \Delta u \nabla
  \bar u \cdot \nabla h+2 \int |\nabla u |^2 \Delta h+2\Re \int \bar
  u\nabla u \nabla \Delta h \\
 & \,\,\,\,\,\,\,{}+2\int \eps |u|^{p-1}
  \nabla(|u|^2) \nabla h +2 \int \eps |u|^{p+1} \Delta h\\
                          & = -4 \Re \int \Delta u \nabla
  \bar u \cdot \nabla h+ 2\int |\nabla u |^2 \Delta h- \int |u|^2
  \Delta^2 h  + \int 2\eps(1-\frac{2 }{p+1}) |u|^{p+1} \Delta h\,.
\end{align*}
Integrating by parts again,
$$
 \int \Delta u \nabla
  \bar u \cdot \nabla h=\int_{\partial \Omega} \nabla \bar u\cdot
 \nabla h \partial_n u-\int \nabla(\nabla \bar u\cdot \nabla h)\cdot
 \nabla u,
$$
and, as $u_{\partial \Omega}=0$ implies $ \partial_\tau u_{\partial
  \Omega}=0$,
\begin{align*}
  2\Re \int \Delta u \nabla
  \bar u \cdot \nabla h & = 2 \int_{\partial \Omega} (\partial_n h)
  |\partial_n u|^2-\int \nabla h\cdot \nabla(|\nabla u|^2)- 2 \int
  \text{Hess} (h) (\nabla u,\nabla \bar u)\\
& =  \int_{\partial \Omega} (\partial_n h)
  |\partial_n u|^2 +\int |\nabla u|^2\Delta h-2  \int
  \text{Hess} (h) (\nabla u,\nabla \bar u)\\
\end{align*}
and finally we obtained
\begin{align*}
   \frac{d^2}{dt^2} M_h(t) & = - \int
   |u|^2 \Delta^2 h  + 2\eps (1-\frac 2 {p+1}) \int |u|^{p+1} \Delta h-2  \int_{\partial \Omega} (\partial_n h)
  |\partial_n u|^2 + 4 \int
  \text{Hess} (h) (\nabla u,\nabla \bar u)
\end{align*}
where we can switch sign for the boundary term if we integrate with
the inner normal of the domain (outer normal of the obstacle !),
retaining the same notation $\partial_n$:
\begin{equation}
  \label{eq:virieldomain}
   \frac{d^2}{dt^2} M_h(t)  = - \int
   |u|^2 \Delta^2 h  + 2 \eps \frac{p-1}{p+1} |u|^{p+1} \Delta h+2  \int_{\partial
 (\R^n\setminus \Omega)} (\partial_n h)
  |\partial_n u|^2 + 4  \int
  \text{Hess} (h) (\nabla u,\nabla \bar u).
\end{equation}
One immediately infers that if $\eps=0,1$, and $h$ is chosen to be the
distance to a point inside the star-shaped obstacle, one controls the
boundary term: let $h(x)=|x|=r$ where the origin $O$ is such that the
obstacle is star-shaped with respect to $O$, then
$$
\Delta h = \frac {n-1} r , \,\,\, \Delta^2 h = \frac{(n-1)(3-n)} {r^3} , \,\,\, \text{Hess} h\geq 0\,\text{ and } n(x)\cdot x>0,
$$
so that $0< \partial_n h\leq 1$. This provides control of the boundary
term in Proposition \ref{p1d}. Now, if one picks
$h(x)=\sqrt{1+|x|^2}$, we still retain $\Delta h\geq 0$ and $\Delta^2
h\leq 0$, and moreover,
$$
\frac{|\nabla u|^2}{h(x)^3}\lesssim  \text{Hess} (h) (\nabla u,\nabla
\bar u),
$$
which implies the local smoothing part in Proposition \ref{p1d}.

 From Proposition \ref{p1d} and Theorem
\ref{t1d}, we immediately deduce Proposition \ref{p2d}.
\begin{rem}
   When $n=2$, the
$\Delta^2 h$ term has the wrong sign: one can only write (without any
attempt to optimize !)
$$
\int_0^T |\int |u|^2 \Delta^2 h |\lesssim C(\Omega) T
\|u_0\|^2_{L^2(\Omega)}.
$$
\end{rem}

Assume now $\eps=0$ and $\Omega$ is domain, either exterior of a
compact set such that there is
no trapped ray, or bounded (boundedness is
not essential). Pick a part of the boundary $P$ where one has local
coordinates such that the normal is a coordinate, and define
$h(x)=d(x,\partial\Omega) \phi(x)$ where $\phi$ is a smooth cut-off to
this coordinate patch, such that on a strip close to the boundary
$\phi$ only depends on the tangential variables. Hence $\partial_n
h\geq 0$ on the boundary part of the patch, and is actually $1$ on a
smaller subset $Q\subset P$, and we control
$$
\int_{Q} |\partial_n u|^2 \leq \int_{P} |\partial_n u|^2.
$$ 
Now, as $M_h$ is controlled by $\|u\|^2_{\dot H^\frac 1 2 }\lesssim
\|u\|_{L^2} \|u\|_{\dot H^1}$, we get using
\eqref{eq:virieldomain}
, with
$S$ a (compact) strip close to the boundary,
$$
\int_0^T \int_{Q} \|\partial_n u\|^2 \lesssim \int_0^T \| u\|^2_{H^1(S)}
ds+\sup_{(0,T)} \|u\|^2_{\dot H^\frac 1 2},
$$
as all the $h$ terms are bounded. Patching together a finite number of
local coordinates patches, we control the entire boundary term,
\begin{equation}
  \label{eq:tracedomainviriel}
  \int_0^T \int_{\partial \Omega} \|\partial_n u\|^2 \lesssim \int_0^T \| u\|^2_{H^1(S)}
ds+\sup_{(0,T)} \|u\|^2_{\dot H^\frac 1 2}.
\end{equation}
Now, on the exterior of a compact set with no trapped rays, local
smoothing holds (\cite{BurqDuke}) and we obtain
a global in time control
\begin{equation}
  \label{eq:tracedomainviriel3D}
  \int_0^T \int_{\partial \Omega} \|\partial_n u\|^2 \lesssim
  \sup_{(0,T)} \|u\|^2_{\dot H^\frac 1 2 }\lesssim
  \|u_0\|_{L^2}\|u_0\|_{\dot H^1}.
\end{equation}
Combining this with Theorem \ref{t1d}, we deduce Theorem \ref{t2d} in
the non star-shaped case, for the linear equation, provided we can
replace the right-hand side $\|u_0\|^2_{L^2_x}\|u_0\|^2_{\dot H^\frac
  1 2}$ by an $\dot H^\frac 1 4$ norm. Let us consider the $\R^n$
case: assume we apply our estimate to
a spectrally localized data $\Delta_j u_0$: then, on the right-hand side,
$$
\|\Delta_j u_0\|^2_{L^2_x}\|\Delta_j u_0\|^2_{\dot H^\frac
  1 2} \sim 2^j \|\Delta_j u_0\|^4_{L^2_x}.
$$
On the left-hand side, we get, summing in $j$,
$$
\sum_j \| |\nabla_x|^{\frac{3-n} 2} (|\Delta_j u|^2)\|_{L^2_{t,x}}
\lesssim \| u_0 \|^2_{\dot H^\frac 1 4}.
$$
As
$$
\| \sum_j |\Delta_j u|^2\|_{\dot H^\frac{3-n} 2}\lesssim \sum_j \| |\nabla_x|^{\frac{3-n} 2} (|\Delta_j u_0|^2)\|_{L^2_{t,x}},
$$
when $n=3$ we are done, using the equivalence of the $L^p$ norm of $u$
with
the $L^p$ norm of its square function, for $p=4$. For $n\neq 3$, one may
decompose $|u|^2$ as a sum of a paraproduct and a reminder: our
previous computation deals with the reminder, while the paraproduct
term can be dealt with by applying the bilinear version of the
estimate to $\Delta_j u$ and $\Delta_k u$ with $k<<j$. We leave the
details to the reader. The case of the exterior domain is dealt with
in a similar way (using the Dirichlet Laplacian spectral localization !).
\begin{rem}
  Here and hereafter we define the fractional Sobolev spaces through
  the spectral localization. They do coincide with the usual ones in
  the range we are interested in, see \cite{Tr83}, and the usual
  properties of the square function extend as well. Alternatively, one
  could define the localization through the heat flow and re-derive all
  required properties by hand, or define all spaces by summing a part
  which is localized close to the boundary (for which one may use all
  the known spectral properties on a bounded domain) and a part which
  is localized away from the boundary (and therefore belongs to the
  usual spaces defined on $\R^n$). This latter approach is essentially
  a poor man's version of the (spatial) localization property of Triebel-Lizorkin spaces, a key point in \cite{Tr83}.
\end{rem}
Now we proceed with the bounded domain $\Omega$: call $v$ the
extension of $u$ by $0$ outside $\Omega$. Then we just proved
\begin{equation}
  \label{eq:rhsdb}
  \| |\nabla|^{\frac {3-n} 2} (|v|^2)\|^2_{L^2(0,T;L^2)} \lesssim \sup_{(0,T)}
\left ( \|u\|^2_2( \|u\|_2 \|u\|_{\dot H^1} +T \|u\|^2_{H^1}) \right).
\end{equation}
Now, assume $u$ to be spectrally localized at (dyadic) $N$:
$u=\phi(N^{-2} \Delta) u$ with $\phi\in C^\infty_0$ and the operator
$\phi(\Delta)$ is defined by functional calculus through the spectral
measure. Now, picking $T$ to be of size $N^{-1}$, the right hand side
in \eqref{eq:rhsdb}  will be bounded by
$\|u\|^4_{H^\frac 1 4}$.
  Consider an interval $[0,1]$, by
subdivision, one gets
$$
\| |\nabla|^{\frac {3-n} 2} (|v|^2)\|^2_{L^2(0,1;L^2)} \lesssim N
\sup_{(0,1)} \|u\|^4_{H^\frac 1 4}\lesssim \|u\|^4_{H^\frac 1 2}.
$$
Now, by \cite{Tr83}, the Sobolev norm of $|u|^2$ is equally the infimum
over all extensions to $\R^n$, hence
$$
\| |\nabla|^\frac {3-n} 2 (|u|^2)\|^2_{L^2(0,1;L^2(\Omega))} \lesssim N
\sup_{(0,1)} \|u\|^4_{H^\frac 1 4}\lesssim \|u_0\|^4_{H^\frac 1 2}.
$$
Finally, one may freely pass from this inequality (which holds for a
spectrally localized function $u$) to the general one with $u_0\in
H^\frac 1 2(\Omega)$ by summing the dyadic pieces (built on the
spectral localization).
\subsection{Scattering in $\R^n$}
Rather than developing the entire theory for all nonlinearities with $1+\frac 4
n<p<1+\frac 4 {n-2}$, we focus on a couple of explicit examples. We
feel that they are generic and provide a straightforward illustration
of required techniques. As on a domain, we set
$$
M=\int |u|^2 \,dx \text{ and } E=\frac 1 2 \int |\nabla u|^2 \,
dx+\frac 1 {p+1} \int |u|^{p+1} \,dx,
$$
which are both conserved quantities.  Consider, for $n=2$, and on $\R^2$,
\begin{equation}
  \label{eq:nls5-2D}
  i\partial_t u+\Delta u=|u|^4 u, \text{ with } u_{|t=0}=u_0\in H^1.
\end{equation}
Local well-posedness can easily be obtained for $\dot H^\frac 1 2$ datum,
and scattering requires control of appropriate space-time norms. An
important feature of the local well-posedness result is that one may
use Sobolev embedding in the course of the proof (as a consequence of
the supercritical exponent with respect to $L^2$). On the other hand,
our a priori bound \eqref{eq:borneapriori}, together with Sobolev
embedding, yields
$$
\|u\|_{L^4_t L^8_x} \lesssim E^\frac 1 8 M^\frac 3 8.
$$
By Gagliardo-Nirenberg inequality, using $u\in L^\infty_t (\dot H^1)$, one gets
$$
\|u\|_{L^6_t L^{12}_x} \lesssim E^\frac 1 4 M^\frac 1 4,
$$
and this quantity scales like the $L^\infty_t(\dot H^\frac 1 2)$ norm,
which is the scaling invariant norm.
Hence, by H\"older and the Leibniz rule, using the $L^6_tL^{12}_x$ norm on four factors
and either $L^\infty_t(L^2_x)$ or $L^\infty_t(\dot H^1)$ norm
on one factor, we get two bounds,
\begin{equation*}
  \| |u|^4 u\|_{L^\frac 3 2_t(L^\frac  6 5_x)} \lesssim E M^\frac 32
  \,\,\text{ and }\,\,\| |u|^4 u\|_{L^\frac 3 2_t(\dot W^1_\frac  6 5)}  \lesssim  E^\frac 32 M.
\end{equation*}
Scattering in $L^2$ and $H^1$ follows immediately by Duhamel, as
$(\frac 3 2,\frac 6 5)$ is a sharp Strichartz admissible pair. By
interpolation, one can obtain scattering for all $H^s$ with $0<s<1$.

One may want to take advantage of the nonlinear part of our a priori
bound \eqref{eq:borneapriori}; in fact, one has
$$
\||u|^4\|_{L^2_t (\dot H^{-\frac 1 2})}\lesssim E^\frac 1 4 M^\frac 3 4.
$$
Combining this with the energy bound, one may prove that
$$
\| |u|^4 u\|_{L^2_t \dot B^{\frac 1 2,2}_1} \lesssim E^{\frac 5 4}
M^\frac 3 4.
$$
If the end-point Strichartz estimate were true with $n=2$, then one
gets a better polynomial bound than the previous one. The lack of the
end point may be routed around to obtain the last bound, with a small
$\epsilon$ loss in the power of $M$.

We now provide another example, which illustrates that one may not
always use the nonlinear part, and that estimating ``in one shot'' the
right space-time norm is not necessarily doable, especially when close
to the $L^2$-critical case. Consider, for $n=1$, and on $\R$,
\begin{equation}
  \label{eq:nls5-1D}
  i\partial_t u+\Delta u=|u|^5 u, \text{ with } u_{|t=0}=u_0\in H^1.
\end{equation}
Local well-posedness can easily be obtained for $\dot H^\frac 3 {10}$ datum.
The linear part of our a priori bound \eqref{eq:borneapriori}, can be interpolated with
the mass conservation, and yields
$$
\||u|^2 \|_{L^3_t L^\infty_x} \lesssim E^\frac 1 {6} M^\frac 5 {6}.
$$
Hence, by product laws, using the $L^6_tL^{\infty}_x$ norm on three
factors, the mass on one and Duhamel, we get (estimating the nonlinear
part in $L^1_t L^2_x$, which may not be the optimal choice !)
$$
\|u\|_{L^4(t_1,t_2; L^\infty_x)} \lesssim M^\frac 1 2+ \|
u\|^3_{L^6(t_1,t_2;L^\infty_x)} \|u\|^2_{L^4(t_1,t_2;L^\infty_x)}.
$$
Assume that $(t_1,t_2)$ is such that
$$
\|u\|^3_{L^6(t_1,t_2;L^\infty_x)} M^\frac 1 2 \lesssim 1/10,
$$
then
$$
\|u\|_{L^4(t_1,t_2; L^\infty_x)} \lesssim 2 M^\frac 1 2.
$$
Splitting the $L^6_t(L^\infty_x)$ in a finite number $N$ of small
increments of equal size $S$, such that 
$$
S^6 M \sim 10^{-2} \text{ and } N
S^6=\|u\|^6_{L^6_t(L^\infty_x)},
$$
 one controls the $L^4_tL^\infty_x$ norm, and the number of increments
 is $N=100 M\|u\|^6_{L^6_t L^\infty_x}$; namely,
$$
\|u\|^4_{L^4_t( L^\infty_x)} \lesssim 16 N M^2\lesssim M^3 (E^\frac 1
{6} M^\frac 5 {6})^3,
$$
and finally
$$
\|u\|_{L^4_t( L^\infty_x)} \lesssim M^\frac {11} 8 E^\frac 1
{8}.
$$
Scattering in $L^2_x$ follows immediately by Duhamel. Now, scattering in
$\dot H^1$ follows by the exact same computation, using Leibniz rule, namely
$$
\|\partial_x u\|_{L^4(t_1,t_2; L^\infty_x)} \lesssim E^{\frac{1} 2}+ \|
u\|^3_{L^6(t_1,t_2;L^\infty_x)}
\|u\|_{L^4(t_1,t_2;L^\infty_x)}\||\partial_x u\|_{L^4(t_1,t_2;L^\infty_x)};
$$
by interpolation, one may then obtain scattering in any $\dot H^s$ for
$0<s<1$. On the other hand, attempts to use the $\int_{t,x}
|u|^9 $ nonlinear a priori bound seem to be doomed by scaling
considerations. For $p\geq 13$, however, it becomes immediately
relevant (notice that for $p=13$, $(p+3)/(p-1)=4/3$).
\begin{rem}
  Informally, for all dimensions $n\geq 1$ one may obtain scattering
  of the full range $1+\frac 4 n<p< 1+\frac 4 {n-2}$ from (the linear
  part of) estimate \eqref{eq:borneapriori}. This can be seen through
  scaling considerations: one is given an a priori space-time bound at the level
  of the $\dot H^\frac 1 4$ norm. Through interpolation with the
  relevant bound (either energy or mass), one retrieves a
  scale-invariant space-time bound. As the equation is
  $L^2$-supercritical, the fixed point argument is using even a tiny
  bit of Sobolev embedding to estimate the nonlinearity, and this is
  enough to insert the a priori estimate and close a true scale
  invariant Strichartz bound.
\end{rem}
\subsection{Nonlinear equation on a domain}
We first deal with Theorem \ref{existence3D}. Notice that the
interesting case is $3\leq p<5$, and we assume for the rest of the
proof that $p$ is close to $5$, which is the most difficult case. Let
us set notations: for any $1\leq q\leq +\infty$, $L^q_t$ denotes a global in
time norm, while $L^q_T$ denotes the norm on a finite time interval
$(0,T)$. Moreover, any implicit constant in a $\lesssim$ sign does not
depend on $T$ (in other words, time dependence is explicitly
tracked). 

We start with linear estimates on the homogeneous and inhomogeneous
equation.
\begin{lemme}
\label{bourbaki}
  Let $S(t)$ denote the linear flow for the Schr\"odinger equation on
  an exterior domain $\Omega$ which satisfies the non trapping
  condition and let $s\geq 0$. Then,
  \begin{equation}
\label{domainl4s1}
  \|S(t)u_0\|_{L^4_t(\dot W^{s,4}_0)} \lesssim \|u_0\|_{\dot
    H^{s+\frac 1 4}_0(\Omega)}.
  \end{equation}
Denote by $w$ the solution of the inhomogeneous equation,
e.g. $w=\int_0^t S(t-s) f(s)\,ds$,
  \begin{equation}
    \label{eq:domaindl4s0}
\|w\|_{C_t(\dot H^{s+\frac 1 2}_0)}+ \|w \|_{L^4_t(\dot W^{s,4})}\lesssim \|f\|_{L^{\frac 4 3}_t \dot
      W^{s+\frac 1 2}_{\frac 4 3}}.
  \end{equation}
Let $\chi_1$ and $\chi_2$ be two
smooth cut-off functions which are such that $\chi_1=1$ on a
ball $B_1$ such that $\R^3\setminus \Omega \subset B_1$, $\chi_1=0$
outside of $2B_1$ and $\chi_2=1$ on $8B_1$, $\chi_2=0$ outside
$9B_1$. Then
\begin{equation}
\label{domaincut}
  \|\chi_2 S(t)u_0 \|_{L^2_t(H^\frac 3 2)}+ \|(1-\chi_1) S(t)u_0\|_{L^4_t(L^\infty_x)} \lesssim \|u_0\|_{H^1_0},
\end{equation}
and
\begin{equation}
\label{domaincutinhom}
\|w \|_{L^4_t(\dot W^{\frac 3 4,4}_0)}+ \|\chi_2 w \|_{L^2_t(H^\frac 3
  2)}+ \|(1-\chi_1) w\|_{L^4_t(L^\infty_x)} \lesssim \|f \|_{L^1_t(H^1_0)}.
\end{equation}
\end{lemme}
Recall that from Theorem
\ref{t2d}, we have an estimate on the linear flow $S(t)$:
\begin{equation}
\label{domainl4s14}
  \|S(t)u_0\|_{L^4_{t,x}} \lesssim \|u_0\|_{\dot H^\frac 1 4(\Omega)}.
\end{equation}
One may shift regularity by $s$ and obtain \eqref{domainl4s1} using fractional powers of the Laplacian and equivalence of norms on
domains (\cite{Tr83}); by the standard $TT^\star$ argument, we also obtain
\eqref{eq:domaindl4s0}, which may again be shifted in regularity
should it be necessary.
\begin{rem}
  Notice that \eqref{domainl4s1} barely fails to provide
control of $L^4_t(L^\infty_x)$. One has to find an appropriate way to
turn around this problem in order to deal with the nonlinear equation. Informally, one may use local smoothing
estimates close to the boundary, and Strichartz estimate for the usual
Laplacian on $\R^3$ away from it. The subcriticality with respect to
$H^1$ of the nonlinear equation will compensate the weakness of the local smoothing estimate.
\end{rem}
Let us now prove \eqref{domaincut}. Notice that multiplying by $\chi_1$ or $\chi_2$ localizes close to
$\partial\Omega$ and $\chi_2=1$ on the support of $\chi_1$. Denote
$u_L=S(t)u_0$. Then, the estimate on $\chi_2 u_L$ in \eqref{domaincut} follows immediately, by local smoothing (see \cite{BGTexterieur}).

Consider now $(1-\chi_1) u_L$: it solves
\begin{equation}
\label{eqchi1}
  i\partial_t (1-\chi_1) u_L+\Delta (1-\chi_1) u_L=[\chi_1,\Delta] u_L,
\end{equation}
and the equation on the left is now set on $\R^3$. We proceed with a useful abstract
lemma which is a simple consequence of a maximal function estimate due
to Christ and Kiselev (\cite{CK}. See also \cite{BPjfa} for a direct
proof without Whitney decompositions).
\begin{lemme}
\label{lemCK}
  Let $U(t)$ be a one parameter group of operators, $1\leq r<
  q\leq+\infty$, $H$ an Hilbert space and $B_r$ and $B_q$ two Banach
  spaces. Suppose that
$$
\|U(t) \varphi\|_{L^q_t(B_q)}\lesssim \|\varphi\|_H \,\,\text{ and }
\,\, \|\int_s U(-s) g \,ds\|_{H}\lesssim \|g \|_{L^r_t(B_r)},
$$
then
\begin{equation}
  \label{eq:CKtype}
  \|\int_{s<t} U(t-s) g(s) \,ds\|_{L^q_t(B_q)}\lesssim \|g\|_{L^{r}(B_r)}.
\end{equation}
\end{lemme}
Now, we pick $U(t)=S(t)$, $L^q_t(B_q)=L^3_t(\dot W^{1,\frac {18} 5})$,
 $L^r_t(B_r)=L^2_t(H^\frac 1 2_{\text{comp}})$ and $H=H^1_0$. Then the
 homogeneous estimate in Lemma \ref{lemCK} is a Strichartz estimate
 (in $\R^3$, with Strichartz pair $(3,\frac{18} 5)$) while the
 inhomogeneous estimate is the dual version of the local smoothing
 (shifted at the right regularity, see again
 \cite{BGTexterieur}). Therefore, applying it to the inhomogeneous
 part of the solution to \eqref{eqchi1}, we get
$$
\|(1-\chi_1)u_L \|_{L^3_t( \dot W^{1,\frac {18}5 })}\lesssim
\|u_0\|_{H^1_0}+\|[\chi_1,\Delta] u_L\|_{L^2_t(\dot H^\frac 1 2)}.
$$
 For a given
function $\phi$, Gagliardo-Nirenberg inequality reads
$$ 
\|
\phi\|_{L^\infty_x} \lesssim \|\phi\|^\frac 1 4_{L^6_x}
\|\phi\|^\frac 3 4_{\dot W^{1,\frac {18}5 }},
$$
which yields
$$ 
\|(1-\chi_1) u_L\|_{L^4_t L^\infty_x} \lesssim \|(1-\chi_1)u_L\|^\frac
1 4_{L^\infty_t \dot H^1_0}
\|(1-\chi_1)u_L \|^\frac 3 4_{L^3_t \dot W^{1,\frac {18}5 }},
$$
and
\begin{equation}
  \|(1-\chi_1) u_L\|_{L^4_t(L^\infty_x)}\lesssim
  \|u_0\|_{H^1_0}+\|[\chi_1,\Delta] u_L\|_{L^2_t(\dot H^\frac 1 2)}.
\end{equation}
Finally, as $\|[\chi_1,\Delta] u_L\|_{L^2_t(\dot H^\frac 1 2)}\lesssim
\|\chi_2 u_L\|_{L^2_t(\dot H^\frac 3 2)}$, we have obtained
\begin{equation}
  \|(1-\chi_1) u_L\|_{L^4_t(L^\infty_x)}\lesssim
  \|u_0\|_{H^1_0}.
\end{equation}
Consider the inhomogeneous equation, $i\partial w+\Delta w=f$, with
$w_{|t=0}=0$. Assume $f\in L^1_T(H^1_0)$, then by using local smoothing and
our Strichartz estimate \eqref{domainl4s1} on $S(t)$ (with $s=\frac 3 4$), and the Duhamel
representation of $w$, we get
\begin{equation}
  \|\chi_2 w\|_{L^2_T( H^\frac 3 2)}+\|w\|_{L^4_T( W^{\frac 3
    4,4})}\lesssim \|f\|_{L^1_T(H^1_0)}.
\end{equation}
Again, consider $(1-\chi_1) w$, solution to
\begin{equation}
    i\partial_t (1-\chi_1) w+\Delta (1-\chi_1) w=[\chi_1,\Delta] w+(1-\chi_1)f,
\end{equation}
exactly as before we get
\begin{equation}
  \|(1-\chi_1) w\|_{L^4_T(L^\infty_x)}\lesssim
  \|f\|_{L^1_T(H^1_0)},
\end{equation}
which ends the proof of Lemma \ref{bourbaki}.

We are now ready to set up a fixed point procedure for equation
\eqref{equd} in the Banach space
\begin{equation}
  \label{eq:X}
  X=\{ u\tq u\in C_T(H^1_0)\cap L^4_T(W^{\frac 3 4,4}),\,\,\,\chi_2
  u\in L^2_T(H^\frac 3 2),\,\,\,(1-\chi_1) u\in L^4_T(L^\infty_x) \}.
\end{equation}
Such a fixed point is standard (and will be omitted) once the following lemma is proven.
\begin{lemme}
\label{contract}
  Let $f=|u|^{p-1} u-|v|^{p-1} v$ with $u,v\in X$. Then, for $p<5$,
\begin{equation}
  \label{eq:nlmap}
  \|f\|_{L^1_T(H^1_0)}\lesssim T^{0^+} \|u-v\|_X (\|u\|^{p-1}_X+\|v\|^{p-1}_X).
\end{equation}
\end{lemme}
Introduce $\chi_3=1$ on $4B_1$, $\chi_3=0$ outside
$5B_1.$
 Let us start with $\chi_3 f$: due to the support conditions, one may
replace $u$ and $v$ by $\chi_2 u$ and $\chi_2 v$ for as many factors
as we wish. By interpolation between $L^\infty_T(H^1)$ and
$L^2_T(H^\frac 3 2)$ for $\chi_2 u$ and interpolation between
$L^\infty_T(H^1)$ and $L^4_T(W^{\frac 3 4,4})$ for $u$, we have
\begin{equation}
  \label{eq:interX}
  \|\chi_2 u\|_{ L^m_T(H^{1+\frac 1 m})}+\|u\|_{ L^q_T(L^{3r})}
  \lesssim \|u\|^{\ }_X,
\end{equation}
where $\frac 1 r=\frac 1 2-\frac 2 q$, and $m$ (resp. $q$) is to be thought of as
very large (resp. slightly larger than $4$). We proceed using
$H^{1+\frac 1 m} \hookrightarrow W^{1,\lambda}$ with
$1/\lambda=1/2-1/(3m)$ and evaluate $\nabla (\chi_3 f)$: by chain rule, we are
left with $\nabla (\chi_2 g) (\chi_2 h)^{p-1}$ where $g,h$ may be $u$,
$v$ or $u-v$ (with one factor of $u-v$ in the $p$ factors). By H\"older, we obtain
\begin{equation*}
  \|\nabla \chi_3 f\|_{L^\rho_T(L^2_x)}\lesssim \|\chi_2
  g\|_{L^m_T(W^{1,\lambda})} \|\chi_2 h\|^{p-1}_{L^q_T(L^{3r}_x)},
\end{equation*}
with 
$$
\frac 1 \rho =\frac 1 m+\frac{p-1} q \text{ and } \frac 1 m=\frac
{p-1}{r}.
$$
Let $p-1=4-\varepsilon$, and pick $m$ such that $\varepsilon>2/m$,
then $\rho>1$ and we recover the correct mapping, with a factor
$T^{1-\frac 1 \rho}$ and $1-\frac 1 \rho=\frac{\varepsilon}4-\frac 1 {2m}$:
\begin{equation*}
  \|\chi_3 f\|_{L^1_T(H^1_0)}\lesssim T^{1-\frac 1 \rho} \|\chi_2
  (u-v)\|_{X} (\|\chi_2 u\|^{p-1}_{X}+\|\chi_2 v\|^{p-1}_{X}).
\end{equation*}
We now turn to $(1-\chi_3)f$; exactly as before, one may consider
$(1-\chi_1) u$ and $(1-\chi_1)v$ rather than $u$ and $v$. Then one has
trivially
\begin{equation*}
  \|(1-\chi_3) f\|_{L^1_T(H^1_0)}\lesssim T^{\frac \varepsilon 4} \|(1-\chi_1)  (u-v)\|_{X} (\|(1-\chi_1) u\|^{p-1}_{X}+\|(1-\chi_1) v\|^{p-1}_{X}),
\end{equation*}
using $u,v\in L^\infty_T(H^1_0)$ on one factor and $(1-\chi_1)(u,v)\in
L^4_T(L^\infty_x)$ on the $4-\varepsilon$ remaining factors. This
achieves the proof of Lemma \ref{contract}.

Local existence and uniqueness in $X$ follows by standard
arguments. Moreover, the local time of existence $T$ is such that
\begin{equation}
  \label{eq:localT}
  T^{\frac \varepsilon 4-\frac 1 {2m}}
  \|u_0\|^{4-\varepsilon}_{H^1_0}\lesssim 1,
\end{equation}
and one may use the conservation of energy to obtain global existence
in the defocusing case. This achieves the proof of Theorem \ref{existence3D}.

We turn to the scattering problem. As in the previous section, we only
provide an explicit example rather than
the best possible general case, for the sake of the exposition. We consider
the defocusing cubic equation on a $3D$ exterior of a domain $\Omega$,
\begin{equation}
\label{eq3-3d}
  i\partial_t u+ \Delta u= |u|^{2} u, \text{ with }
  u_{|\partial \Omega} =0,\,\,\, u_{t=0}=u_0\in H^1_0(\Omega),
\end{equation}
and require $\Omega$ to be star-shaped. Let
$\chi_1,\chi_2,\chi_3 $ be smooth cut-off functions close to the boundary
$\partial \Omega$, with $\chi_1\leq \chi_3\leq \chi_2$ as in the
existence proof we just completed.

Recall that we have two different nonlinear estimates which are valid
for all times.
\begin{itemize}
\item From Proposition \ref{p1d}, we control a local smoothing type
  quantity at the level of $H^\frac 1 2$ regularity on the data: 
\begin{equation}
  \label{eq:localnonlinear}
  \int_0^{+\infty} \|\chi_2 u\|^2_{\dot H^1_0}\,dt \lesssim
  M(u_0)^\frac 1 2 E^\frac 1 2(u_0).
\end{equation}
\item From Proposition \ref{p2d}, we control a Strichartz-like norm,
  at the level of regularity $H^\frac 1 4$ on the data:
  \begin{equation}
    \label{eq:L4scat}
    \| u\|_{L^4_t(L^4(\Omega)} \lesssim M^\frac 3 8 E^\frac 1 8.
  \end{equation}
\end{itemize}
Ultimately, we aim at controlling space-time norms at the level of
$H^1$ regularity on the data. We start by bootstrapping our relatively
weak control \eqref{eq:L4scat} into a somewhat stronger estimate at
the level of $H^\frac 1 2$ regularity.
\begin{lemme}
\label{lemboot1}
  Let $u$ be a solution of \eqref{eq3-3d}. Then
  \begin{equation}
    \label{eq:bootstrap1}
     \chi_1 u \in L^4_t (\dot W^{\frac 1 4,4}) \,\,\text{ and }\,\,
(1-\chi_1) u \in L^3_t(\dot W^{\frac 1 2,\frac{18} 5})\cap L^\frac
{12} 5_t(\dot W^{\frac 1 2,\frac 9 2}).
  \end{equation}
As a consequence, the solution $u$ scatters in $H^\frac 1 2$.
\end{lemme}
In order to prove the Lemma, we again split the equation, treating
differently the neighborhood of the boundary (where local smoothing is
most efficient) and spatial infinity (where Strichartz estimates for
the free propagator are available). Consider $\chi_1 u$, which is a solution to
\begin{equation}
  i\partial_t \chi_1 u+ \Delta \chi_1 u=\chi_1 |u|^{2} u -[\chi_1,\Delta] u=f.
\end{equation}
On the nonlinear part, which is compactly supported, we use
\eqref{eq:localnonlinear} on one factor, and $L^\infty_t(H^1_0)$ on
the other two, while the commutator term is easily controlled by
$\|\chi_2 u\|_{L^2(\dot H^1_0)}$; hence, globally in time,
$$
\|f\|_{L^2_t(L^2_{\text{comp}})} \lesssim \|\chi_2 u\|_{L^2_t(\dot
  H^1_0)} (E(u_0)+1).
$$
Then, we apply Lemma \ref{lemCK}, this time with $H=\dot H^\frac 1
2$, $L^q_t(B_q)=L^4_t(\dot W^{\frac 1 4,4})$ (this from interpolation
between \eqref{domainl4s1} and \eqref{domainl4s14})) and
$ L^p_t(B_p)=L^2_t(L^2_{\text{comp}})$ (which is dual local smoothing at regularity
$H^\frac 1 2$); we obtain $\chi_1 u \in L^4_t (\dot W^{\frac 1 4,4})$.

Let us deal with $(1-\chi_1) u$, which is solution to
\begin{equation}
\label{eq:extg2}
  i\partial_t (1-\chi_1) u+ \Delta (1-\chi_1) u=(1-\chi_1) |u|^{2} u +[\chi_1,\Delta] u.
\end{equation}
Recall that we are now on the whole space $\R^3$. The commutator term
is dealt with exactly as with the previous part, and is therefore
$L^2_t(L^2_{\text{comp}})$. Another application of Lemma \ref{lemCK}
but with $L^q_t(B_q)$ being either $L^3_t(\dot W^{\frac 1 2,\frac{18}
  5})$ or $L^\frac {12} 5_t(\dot W^{\frac 1 2,\frac 9 2})$ (both $(3,\frac {18}
5)$ and $(\frac {12} 5,\frac 9 2)$ are Strichartz pairs for the free space)
yields the claim for the Duhamel term coming from the commutator. One
the other hand,
$$
\| (1-\chi_1)|u|^2 u\|_{L^2_t(\dot W^{1,1})} \lesssim
\|u\|_{L^\infty_t(\dot H^1)} \|u\|^2_{L^4_{t,x}} \lesssim M^{\frac 3
  4} E^{\frac 3 4},
$$
where we used our other a priori control \eqref{eq:L4scat} and the
energy conservation. From the embedding $L^2_t(\dot W^{1,1})\hookrightarrow L^2_t(\dot W^{\frac 1
  2,\frac 6 5})$, we may apply the dual end-point Strichartz estimate on the
nonlinear term and finally $(1-\chi_1) u \in L^3_t(\dot W^{\frac 1
  2,\frac{18} 5})\cap L^\frac {12} 5_t(\dot W^{\frac 1 2,\frac 9 2})$ which
achieves the proof of \eqref{eq:bootstrap1}.

From $(1-\chi_1)u \in L^3_t(\dot W^{\frac 1
  2,\frac{18}5})\hookrightarrow L^3(L^9_x)$, we infer by Leibniz,
H\"older ($\frac 1 2=\frac 2 9+\frac 5 {18}$) and interpolation,
\begin{equation}
  \label{eq:scatext}
  ((1-\chi_1) u)^3 \in L^1_t(\dot H^\frac 1 2).
\end{equation}
Going back to the equation on $u$ and splitting the source term $|u|^2
u$ as $g_1=\chi_2 |u|^2 u$ and $g_2=(1-\chi_2) |u|^2 u$, we have for
the same reason as before (and globally in time !) $g_1\in
L^2_t(L^2_{\text{comp}})$, while from \eqref{eq:scatext},
$g_2 \in L^1_t(\dot H^{\frac 1 2})$. Then, with $S(t)$ being the
Schr\"odinger group on our domain,
\begin{equation}
  \label{eq:scat12}
  S(-t)u=u_0+\int_0^t S(-s) (g_1+g_2) ds=u_0+\int_0^{+\infty} S(-s) (g_1+g_2) ds-\int_t^{+\infty} S(-s) (g_1+g_2) ds,
\end{equation}
from which scattering in $\dot H^\frac 1 2$ follows: both integral
terms are well-defined in $\dot H^\frac 1 2$, and the second one
vanishes when $t\rightarrow+\infty$. This proves Lemma \ref{lemboot1}.

We aim at bootstrapping this information up to $H^1$ scattering in two
steps. First, we improve our new space-time controls to the level of
$\dot H^\frac 3 4$ regularity: once again, the important point is
to use only global in time bounds.
\begin{lemme}
\label{lemboot2}
  Let $u$ be a solution of \eqref{eq3-3d}. Then
  \begin{equation}
    \label{eq:bootstrap2}
     \chi_2 u \in L^2_t (\dot H^{\frac 5 4}) \,\,\text{ and }\,\,
 u \in L^4_t(\dot W^{\frac 1 2,4}).
  \end{equation}
\end{lemme}
The proof of Lemma \ref{lemboot2} is more delicate than the previous
one, as a splitting time argument (like we did for scattering in
$\R^n$) is required. From the local existence theory, we can easily
get that the $L^4_T(\dot W^{\frac 1 2,4})$ is finite for
$T<+\infty$. We will prove that $T=+\infty$, by using the equation and Duhamel:
\begin{itemize}
\item we start with $g_1=\chi_3 |u|^2 u$: interpolating between $u\in
  L^\infty_t(\dot H^1)$ and $\chi_2 u \in L^2_t(\dot H^1)$, we have
  $\chi_2 u\in L^4_t(\dot H^1)$. On the other hand, from Lemma
  \ref{lemboot1}, $u\in L^4_t(\dot W^{\frac 1 4,4})$; by interpolation
  with $u\in L^4_T(\dot W^{\frac 1 2,4})$ and Sobolev, we get $u^2 \in L^2_T(L^4_x)$ and
$$
\|\chi_3 |u|^2 u\|_{L^\frac 4 3_T(\dot W^{1,\frac 4 3})} \lesssim
\|u\|^\frac 1 2_{L^\infty_t(\dot H^1)}\|\chi_2 u\|^\frac 1 2_{L^2_t(\dot
  H^1)} \|u\|_{L^4_t(\dot W^{\frac 1 4,4})}\|u\|_{L^4_T(\dot W^{\frac
    1 2 ,4})}\,;
$$
\item let us deal with $g_2=(1-\chi_3) |u|^2 u$: interpolating between $u\in
  L^\infty_t(\dot H^1)$ and $(1-\chi_1) u\in L^{3}_t(\dot W^{\frac 1
    2,\frac {18} 5})$ (which we got from Lemma \ref{lemboot1})), we obtain
  $(1-\chi_1) u \in L^6_t(\dot W^{\frac 3 4,\frac {18} 7})$. Recall as well that
Lemma \ref{lemboot1} provides $(1-\chi_1) u\in L^{\frac {12} 5}_t(\dot W^{\frac 1
    2,\frac 9 2})\hookrightarrow L^\frac{12} 5_t(L^{18}_x) $; using this
  information on two factors and the interpolation bound on the third
  one, we get
  \begin{equation}
    \label{eq:L1ext34}
    \|(1-\chi_3 ) |u|^2 u\|_{L^1_t(\dot H^{\frac 3 4})} \lesssim
\|u\|^{\frac 1 2}_{L^\infty_t(\dot H^1)}\|(1-\chi_1) u\|^{\frac 1 2}_{L^{3}_t(\dot
  W^{\frac 1 2,\frac{18}5 })} \|(1-\chi_1)u\|^2_{L^\frac{12} 5_t(\dot W^{\frac 1 2,\frac {9} 2})}.
  \end{equation}
\end{itemize}
Using the equation, Duhamel and \eqref{eq:domaindl4s0} at regularity $s=\frac 1 2$, we have
$$
\|u\|_{L^4_T(\dot W^{\frac 1 2 ,4})} \lesssim \left[\|u_0\|_{\dot H^\frac 3
  4}+
    \|(1-\chi_3 ) |u|^2 u\|_{L^1_t(\dot H^{\frac 3 4})}\right]+ 
M^\frac 1 2 E \|u\|_{L^4_t(\dot W^{\frac 1 4,4})}\|u\|_{L^4_T(\dot W^{\frac
    1 2 ,4})}\,;
$$
the bracket term is finite by \eqref{eq:L1ext34}, and a splitting time
argument performed on the $L^4_t(\dot W^{\frac 1 4,4})$ norm which is
finite provides global in time control of $u\in
L^4_{t}(\dot W^{\frac 1 2,4})$. Using Duhamel, again, on $g_1$ and $g_2$, we also
obtain $\chi_3 u \in L^2_t(\dot H^\frac 5 4)$, globally in time (note
that for $g_1$ we have to resort again to Lemma \ref{lemCK},
combining \eqref{eq:domaindl4s0} and local smoothing). This achieves
the proof of Lemma \ref{lemboot2}. We finally need one last step.
\begin{lemme}
\label{lemboot3}
  Let $u$ be a solution of \eqref{eq3-3d}. Then
  \begin{equation}
    \label{eq:bootstrap3}
     \chi_3 |u|^2u \in L^\frac 4 3_t (\dot W^{\frac 5 4,\frac 4 3}) \,\,\text{ and }\,\,
 (1-\chi_3)|u|^2 u \in L^1_t(\dot H^1_0).
  \end{equation}
\end{lemme}
Again, we proceed differently close to or far from the boundary.
\begin{itemize}
\item On $g_1$, we use $\chi_2 u \in L^2_t(\dot H^\frac 5 4)$ from
  Lemma \ref{lemboot2} and $u\in
  L^8_{t,x}$ (which follows from $u\in L^\infty_t(\dot H^1)$ and $u\in
   L^4_t(\dot W^{\frac 1 2,4})$, again from Lemma \ref{lemboot2}), and obtain
$$
\|g_1\|_{L^\frac 4 3_t(\dot
  W^{\frac 5 4,\frac 4 3})} \lesssim \|\chi_2 u\|_{L^2_t(\dot H^\frac
  5 4)} \|u\|^2_{L^8_{t,x}}.
$$
\item For $g_2$, we need $(1-\chi_3)u\in L^2_t(L^\infty_x)$ which does
  not follow from the Strichartz estimates we already obtained on
  $(1-\chi_1)u$ (missing end-point, not to mention a log). We use \eqref{eq:extg2}, but with the cut
  $\chi_3$ instead of  $\chi_1$:
\begin{equation}
\label{eq:extg3}
  i\partial_t (1-\chi_3) u+ \Delta (1-\chi_3) u=(1-\chi_3) |u|^{2} u +[\chi_3,\Delta] u.
\end{equation}
 For the nonlinear
  part,
$$
\|(1-\chi_3)|u|^2 u\|_{L^2_t(L^\frac 6 5_x)} \lesssim \|u\|_{L^\infty(\dot H^\frac 1
  2)}\|u\|^2_{L^4_{t,x}}.
$$
The commutator term is $L^2_t(L^2_{\text{comp}})$ hence $L^2_t(L^\frac
6 5_x)$ and by Duhamel (with
  the free propagator !) we get by Strichartz $(1-\chi_3)u \in L^2_t(L^6_x)$. On the other hand,
  using again  \eqref{eq:extg3}, $\chi_3 u \in L^2_t(\dot H^\frac 5
  4)$ for the commutator term, and 
  \begin{equation}
    \|(1-\chi_3 ) |u|^2 u\|_{L^1_t(\dot H^{\frac 3 4})} \lesssim
\|u\|^{\frac 1 2}_{L^\infty_t(\dot H^1)}\|(1-\chi_1) u\|^{\frac 1 2}_{L^{3}_t(\dot
  W^{\frac 1 2,\frac{18}5 })} \|(1-\chi_1)u\|^2_{L^\frac{12} 5_t(\dot W^{\frac 1 2,\frac {9} 2})},
  \end{equation}
(this is nothing but \eqref{eq:L1ext34} with $\chi_2$ replaced by
$\chi_3$) we wish to obtain by (free) Strichartz and Duhamel, $(1-\chi_3)u \in
L^2_t(\dot W^{\frac 3 4,6})$. This is indeed the case for the Duhamel
term coming from the nonlinear term. However, one may no longer use
Lemma \ref{lemCK} for the commutator term, and we need in a
crucial way the $L^2_t$ norm. Fortunately enough, we may use
\begin{lemme}[Staffilani-Tataru \cite{ST}]
  Let $x\in \R^n$, $n\geq 3$ and let $f(x,t)$ be compactly supported in space, such
  that $f\in L^2_t(H^{-\frac 1 2})$. Then the solution  $w$ to
  $(i\partial_t +\Delta_x) w=f$ with $w_{|t=0}=0$, is such that
$$
\|w\|_{L^2_t(L^{\frac{2n}{n-2}}_x)} \lesssim \|f\|_{L^2_t(H^{-\frac 1 2})}.
$$
\end{lemme}
In other words, provided the left handside is a local smoothing norm,
one recover the endpoint estimate in addition to the ones provided by
Lemma \ref{lemCK}. A quick inspection of the proof in \cite{ST} allows
one to shift spatial regularity as we need.

Going back to our equation on $(1-\chi_3)u$, we therefore obtain  $(1-\chi_3)u \in
L^2_t(\dot W^{\frac 3 4,6})$; Gagliardo-Nirenberg immediately provides $(1-\chi_3) u\in
  L^2_t(L^\infty_x)$. Combining this with $u\in L^\infty_t(\dot H^1)$,
  we finally get $(1-\chi_2) |u|^2 u\in L^1_t(\dot H^1)$. This
  achieves the proof of Lemma \ref{lemboot3}.
\end{itemize}
From the informations on $g_1$ and $g_2$ provided by Lemma \ref{lemboot3}, we may go back to
\eqref{eq:scat12} and obtain scattering in $ H^1$ like we did in
$H^\frac 1 2$. This achieves the proof of Theorem \ref{scattering3D}.
\begin{rem}
  If one picks $p>3$, scattering in a negative regularity Sobolev
  space may easily be obtained. Bootstrapping appears to be more
  difficult, the numerology working in the wrong
  direction when $p$ gets closer to $5$, as both a priori bounds (the
  smoothing and the $L^4_{t,x}$ are subcritical with respect to scaling).
\end{rem}


\end{document}